\documentclass[reqno]{amsart}
\usepackage{amssymb}
\usepackage[mathscr]{eucal}

\newcommand{\nc}[2]{\newcommand{#1}{#2}}
\newcommand{\rnc}[2]{\renewcommand{#1}{#2}}

\nc{\carh}{\mathbf{h}}
\nc{\kf}{\mathfrak{B}}
\nc{\elem}{\mathbf{e}}
\nc{\ot}{\otimes}
\nc{\tot}{\hat\ot}
\nc{\opl}{\oplus}
\nc{\liea}{\mathfrak{g}}
\nc{\sliea}{\mathfrak{h}}
\nc{\uea}{U(\liea)}
\nc{\quea}{U_{h}(\liea)}
\nc{\sln}{\mathfrak{sl}_{n}(\CC)}
\nc{\slt}{\mathfrak{sl}_{2}(\CC)}
\nc{\slth}{\mathfrak{sl}_{3}(\CC)}
\nc{\sut}{\mathfrak{su}(2)}
\nc{\sun}{\mathfrak{su}(l+1)}
\nc{\qusln}{U_{h}(\sln)}
\nc{\sonn}{\mathfrak{so}_{2n+1}(\CC)}
\nc{\sonnr}{\mathfrak{so}(2n+1)}
\nc{\spn}{\mathfrak{sp}_{2n}(\CC)}
\nc{\spnr}{\mathfrak{sp}(2n)}
\nc{\son}{\mathfrak{so}_{2n}(\CC)}
\nc{\sonr}{\mathfrak{so}(2n)}
\nc{\glv}{\mathfrak{gl}(V)}
\nc{\gll}{\mathfrak{gl}(\mathfrak{g})}
\nc{\gln}{\mathfrak{gl}_{n}(\CC)}
\nc{\glt}{\mathfrak{gl}_{2}(\CC)}
\nc{\slfg}{\mathfrak{sl}_{2N-1}(\CC)}
\nc{\glf}{\mathfrak{gl}_{4}(\CC)}
\nc{\glnr}{\mathfrak{gl}(l+1)}
\nc{\CC}{\mathbb{C}}
\rnc{\AA}{\mathbb{A}}
\nc{\DD}{\mathbb{D}}
\nc{\QQ}{\mathbb{Q}}
\nc{\RR}{\mathbb{R}}
\nc{\ZZ}{\mathbb{Z}}
\nc{\NN}{\mathbb{N}}
\nc{\QQbar}{\overline{\QQ}}
\nc{\ur}{\mathcal{R}}
\nc{\uf}{\mathcal{F}}
\nc{\D}{\Delta}
\rnc{\d}{\delta}
\rnc{\b}{\beta}
\nc{\s}{\sigma}
\rnc{\l}{\lambda}
\rnc{\L}{\Lambda}
\nc{\ep}{\epsilon}
\nc{\fa}{\forall}
\nc{\lah}{\varphi}
\nc{\lgh}{\phi}
\nc{\mapto}{\rightarrow}
\nc{\act}{\triangleright}
\nc{\adact}{\stackrel{\mathrm{ad}}{\act}}
\nc{\radact}{\stackrel{\mathrm{ad}}{\ract}}
\nc{\ract}{\triangleleft}
\nc{\ch}{\mathfrak{h}}
\rnc{\a}{\alpha}
\rnc{\k}{\kappa}
\nc{\x}{\chi}
\nc{\g}{\gamma}
\nc{\G}{\Gamma}
\rnc{\r}{\rho}
\nc{\w}{\omega}
\rnc{\t}{\theta}
\nc{\n}{\eta}
\nc{\VC}{\mathcal{V}}
\nc{\VQ}{\Upsilon}
\nc{\EC}{\mathcal{E}}
\nc{\EQ}{\Omega}
\nc{\lca}{\D^{L}_{A}}
\nc{\rca}{\D^{R}_{A}}
\nc{\lcat}{{\D'}^{L}_{A}}
\nc{\rcat}{{\D'}^{R}_{A}}
\nc{\lcan}{{\D^{n}}^{L}_{A}}
\nc{\rcan}{{\D^{n}}^{R}_{A}}
\nc{\lcata}{{\D^{\ot}}^{L}_{A}}
\nc{\rcata}{{\D^{\ot}}^{R}_{A}}
\nc{\lcad}{{\tilde{\D}}^{L}_{A}}
\nc{\rcad}{{\tilde{\D}}^{R}_{A}}
\nc{\lcae}{{\D^{\EQ}}^{L}_{A}}
\nc{\rcae}{{\D^{\EQ}}^{R}_{A}}
\nc{\lcaed}{{\D^{\tilde{\EQ}}}^{L}_{A}}
\nc{\rcaed}{{\D^{\tilde{\EQ}}}^{R}_{A}}
\nc{\llra}{\Longleftrightarrow}
\newcommand{\pair}[2]{\langle {#1},{#2} \rangle}

\theoremstyle{plain}
\newtheorem{thm}{Theorem}
\newtheorem{prop}[thm]{Proposition}
\newtheorem{lem}{Lemma}

\theoremstyle{definition}
\newtheorem{defn}[thm]{Definition}
\newtheorem{ex}[thm]{Example}

\theoremstyle{remark}
\newtheorem{remark}[thm]{Remark}

\numberwithin{thm}{section}

\nc{\bth}{\begin{thm}}
\nc{\bprop}{\begin{prop}}
\nc{\ble}{\begin{lem}}
\nc{\bde}{\begin{defn}}
\nc{\bex}{\begin{ex}}
\nc{\bre}{\begin{remark}}

\nc{\ethe}{\end{thm}}
\nc{\eprop}{\end{prop}}
\nc{\ele}{\end{lem}}
\nc{\ede}{\end{defn}}
\nc{\eex}{\end{ex}}
\nc{\ere}{\end{remark}}

\nc{\dif}{\mathrm{d}}
\nc{\udif}{\mathrm{D}}

\DeclareMathOperator{\id}{id}

\DeclareMathOperator{\inv}{inv}

\DeclareMathOperator{\rinv}{rinv}

\DeclareMathOperator{\Tr}{Tr}

\DeclareMathOperator{\Img}{Im}

\title[Bicovariant differential calculi on Jordanian quantum groups]
{Classification of bicovariant differential calculi on the Jordanian 
quantum groups $GL_{h,g}(2)$ and $SL_{h}(2)$ and quantum Lie algebras}

\author{Andrew D. Jacobs and J. F. Cornwell}

\address{Department of Physics and Astronomy, University of St.~Andrews,
North Haugh, St.~Andrews, Fife, KY16 9SS, Scotland}
\address{
Current Address of ADJ:~Department of 
Mathematical Sciences, University of Cincinnati, Cincinnati, Ohio. 
e-mail:~Andrew.Jacobs@math.uc.edu}

\begin{document}

\begin{abstract}
We classify all 4-dimensional first order bicovariant calculi on the 
Jordanian quantum group $GL_{h,g}(2)$ and all 3-dimensional first 
order bicovariant calculi on the Jordanian quantum group $SL_{h}(2)$. 
In both cases we assume that the bicovariant bimodules are generated as 
left modules by the differentials of the quantum group generators. It 
is found that there are 3 1-parameter families of 4-dimensional bicovariant 
first order calculi on $GL_{h,g}(2)$ and that there is a single, unique, 3-dimensional 
bicovariant calculus on $SL_{h}(2)$. This 3-dimensional calculus may be 
obtained through a classical-like reduction from any one of the three 
families of 4-dimensional calculi on $GL_{h,g}(2)$. Details 
of the higher order calculi and also the quantum Lie algebras are 
presented for all calculi. The quantum Lie algebra obtained from the 
bicovariant calculus on $SL_{h}(2)$ is shown to be isomorphic to the 
quantum Lie algebra we obtain as an ad-submodule within the 
Jordanian universal enveloping algebra $U_{h}(\slt)$ and also through 
a consideration of the decomposition of the tensor product of two 
copies of the deformed
adjoint module. We also obtain the quantum Killing form for this 
quantum Lie algebra.
\end{abstract}

\maketitle

\section{Introduction}

The program of noncommutative geometry pioneered by 
Connes~\cite{Connes1,Connes2} is based on fundamental results 
in the field of abstract analysis discovered 
in the first half of this century by Gelfand, Kolmogoroff, Naimark, Stone 
and others (a useful historical overview can be found in I.~Segal's review~\cite{Segal} of 
Connes' book~\cite{Connes1}). In particular, Gelfand and Kolmogoroff showed that for a 
locally compact space, the algebra of continuous functions on the 
space is essentially \emph{equivalent} to the space itself. The algebra of 
continuous functions is of course commutative, and in fact a 
$C^{*}$-algebra. We can then reasonably consider the study 
of \emph{noncommutative} $C^{*}$-algebras as some form of noncommutative geometry.
Thus the essential idea is to express the formalism of classical 
geometry as far as possible in the language of commutative algebra, 
and then use this as the paradigm for generalising to the noncommutative 
setting. 

An implementation of this program has been developed 
by Dubois-Violette and co-workers (see the book by Madore~\cite{Madore}, 
and the references 
therein). They generalise an elegant algebraic approach to the 
differential geometry of a smooth manifold introduced by
Koszul~\cite{Koszul} to the case where the commutative algebra
of smooth functions on the manifold is replaced by the noncommutative 
algebra of matrices over some field.

In another direction, quantum groups provide natural candidates for noncommutative 
generalisations of the algebras of smooth functions on 
classical compact Lie groups. Classically 
the algebra of \emph{representative functions} is a dense subalgebra of 
the algebra of all smooth complex valued functions on the group and 
carries the structure of 
a Hopf $\star$-algebra. Dropping the $\star$-structure we obtain the 
coordinate rings of the corresponding complex Lie groups, $SL_{n}(\CC)$, $SO_{n}(\CC)$ and 
$Sp_{2n}(\CC)$. These coordinate rings are generated 
already by the matrix elements of the defining representations of 
these groups (for $GL_{n}(\CC)$ we should adjoin 
$(\det(\mathbf{t}))^{-1}$ where $\mathbf{t}$ is the matrix of matrix 
element functions). Corresponding 
to these classical groups are the well known FRT Hopf algebras 
$A(R)$, introduced by 
Faddeev, Reshetikhin and Takhtajan in~\cite{FRT}. They are 
`quantisations' of the classical coordinate rings. It is these structures upon 
which we should try to develop some sort of noncommutative 
Lie group geometry. The 
pioneering works here are those of Woronowicz~\cite{Wor1,Wor}. In 
particular, in~\cite{Wor} Woronowicz set out a formalism in terms of 
\emph{bicovariant bimodules} which has been 
studied intensively by very many authors since. Let us note that the 
classical differential calculus on Lie groups is bicovariant.

There is a well known `problem' with the bicovariant calculi 
associated with the standard FRT quantum groups other than $GL_{q}(n)$: Their 
dimensions do not agree with the corresponding classical calculi. In 
the particular case of $SL_{q}(n)$ the bicovariant calculus of 
Woronowicz is $n^{2}$-dimensional while the classical calculus is of 
dimension $n^{2}-1$ (the dimension being the dimension of the vector 
space of left-invariant 1-forms). This problem has stimulated some 
authors to consider alternative approaches to the development of 
differential geometry on quantum groups. For example, in~\cite{SS1} Schm\"{u}dgen 
and Sch\"{u}ler consider left-covariant bimodules (developing  
Woronowicz's original approach~\cite{Wor1}) on $SL_{q}(n)$. 
They obtain first order calculi with the classical dimension, however 
for $N\geq 4$ the higher order calculi do not have the correct 
dimension. Another interesting approach was initiated by Faddeev and 
Pyatov\cite{FadPyt}. 
They considered, for the particular case of $SL_{q}(n)$, 
the consequences of relaxing the condition of the classical Leibniz rule 
which is present in 
the bicovariant Woronowicz approach. 
They obtained bicovariant calculi of the correct classical 
dimension at all orders. However, subsequent work by Arutyunov, Isaev 
and Popowicz\cite{AruIsaPop}
suggests that a similar approach cannot be employed for 
the other simple quantum groups $SO_{q}(n)$ and $Sp_{q}(n)$. 

In this paper we consider the original Woronowicz bicovariant 
calculus, but we examine such calculi on the \emph{non-standard} 
quantum groups $GL_{g,h}(2)$ and $SL_{h}(2)$ --- the so-called 
Jordanian quantum groups. In~\cite{Kar1,Kar2}, Karimpour initiated the 
study of bicovariant calculi associated with $SL_{h}(2)$. Here we 
perform a complete classification of all first order bicovariant 
calculi on the quantum groups $GL_{h,g}(2)$ and $SL_{h}(2)$. 
Furthermore we consider the higher order calculi and the 
corresponding quantum Lie algebras. Let us summarise our main 
classification results:

\begin{itemize}
\item
There are three 1-parameter families of 4-dimensional first order 
bicovariant differential calculi on $GL_{h,g}(2)$ whose bimodules of 
forms are generated as left $GL_{h,g}(2)$-modules by the 
differentials of the quantum group generators.
\item
For one value of the parameter, the calculi in the three families are 
the same. This parameter value coincides with the value required for 
a `classical-like' reduction to a 3-dimensional first order 
bicovariant calculus 
on $SL_{h}(2)$ which is shown to be unique.
\item 
For all the calculi the relations in the exterior algebra are 
obtained and are shown to lead to exterior calculi whose dimension 
is classical at all orders.
\item
For all the calculi the relations in the enveloping algebra of the 
quantum Lie algebra are obtained and are shown to lead to PBW-type 
bases.
\end{itemize}

Classically, the Lie algebra of a Lie group is obtained as the vector 
space of tangent vectors at the identity equipped with a Lie bracket defined 
in terms of the left-invariant vector fields on the group manifold. 
The formalism of Woronowicz's bicovariant calculus has a natural construction 
for a `quantum Lie 
algebra' 
generalising the classical construction to the abstract Hopf 
algebra setting. However, for all standard quantum groups 
the quantum Lie algebras so obtained have the `wrong' dimension. This 
has prompted authors such as Sudbery and Delius to look for alternative constructions 
for quantum Lie algebras~\cite{LyuSud,Del1,Del2}. There are two different 
approaches described in their work. Recall that classically the Lie 
algebra $\liea$ is an ad-submodule of the classical adjoint 
$U(\liea)$-module, $U(\liea)$,
and its Lie bracket is the restriction of this classical adjoint 
action to $\liea$. This motivates the first 
approach~\cite{LyuSud,Del1}, \textbf{Approach~1}, in which we look for an
ad-submodule \emph{within} the quantised universal enveloping algebra 
$U_{q}(\liea)$ which has the 
correct dimension and upon which the restriction of the adjoint 
action of $U_{q}(\liea)$ closes.
Lyubashenko and Sudbery employed a result of Joseph and Letzter and 
obtained a quantum Lie algebra in the cases $U_{q}(\sln)$ such that the
coproduct of $U_{q}(\sln)$ applied to the basis of their quantum Lie 
algebra is of a particularly neat form.  The other 
approach\cite{Del2}, \textbf{Approach~2}, (see also the paper by Bremner~\cite{Brem} and the 
earlier paper of
Donin and Gurevich~\cite{DonGur} where the idea appeared originally) 
constructs a quantum Lie algebra independent of any 
embedding into $U_{q}(\liea)$. The idea here is to recall that 
classically the Lie bracket is an intertwiner,
$[,]:\mathrm{ad}\ot\mathrm{ad}\mapto\mathrm{ad}$,
where $\mathrm{ad}$ is the usual adjoint representation of 
$U(\liea)$. The quantum Lie bracket is then obtained as the 
intertwiner of the corresponding $U_{q}(\liea)$-modules. Furthermore,
classically the Killing form is an intertwiner,
$\kf:\mathrm{ad}\ot\mathrm{ad}\mapto\CC$,
and
considering the intertwiner between the tensor product of the two 
(quantum) adjoint modules and the trivial representation an analog of 
the Killing form for the quantum Lie algebras is obtained in a rather 
straightforward manner. The two approaches lead to isomorphic quantum 
Lie algebras. However each has advantages over the other: 
\textbf{Approach 1} 
allows us to see explicitly the relationship between the quantum Lie 
algebra and the quantised universal enveloping algebra and in 
principle allows us to determine the coproduct of $U_{q}(\liea)$ on the quantum 
Lie algebra; \textbf{Approach 2} gives a reasonably simple 
prescription for constructing quantum Lie algebras based on 
computation of inverse Clebsch-Gordan coefficients~\cite{Del3} and also 
a direct means of establishing 
the quantum Killing form.

In the last part of the present paper we pursue this line of 
enquiry starting with the Jordanian quantised enveloping algebra 
$U_{h}(\slt)$. The following results are obtained:

\begin{itemize}
\item
A Jordanian quantum Lie algebra is obtained according to \textbf{Approach 1}  
and the expressions for the coproduct of $U_{h}(\slt)$ 
on the basis elements of the quantum Lie algebra are obtained.
\item
Following \textbf{Approach 2} we obtain a quantum Killing form for 
the quantum Lie algebra.
\item
We observe that the quantum Lie algebras obtained equivalently by 
\textbf{Approach 1} and \textbf{Approach 2} are isomorphic to the 
quantum Lie algebra obtained through the unique bicovariant calculus 
on the Jordanian quantum group $SL_{h}(2)$.
\end{itemize}

The paper is self-contained and organised as follows. In Section~\ref{Jorqugro} we recall 
pertinent definitions and results concerning the Jordanian quantum 
groups $GL_{h,g}(2)$ and $SL_{h}(2)$. In Section~\ref{WorThr} 
Woronowicz's theory is reviewed in a style intended to enlighten the 
details of the classification procedure, due to 
M\"{u}ller-Hoissen~\cite{MH1}, which we present in 
Section~\ref{MHProc}. The classification results appear in 
Sections~\ref{ClasResults1} and~\ref{ClasResults2}. The Jordanian quantised universal 
enveloping algebra is 
recalled in Section~\ref{defuhslt} with the corresponding quantum Lie 
algebras obtained in Sections~\ref{Approach1} and~\ref{Approach2} through 
\textbf{Approach~1} and \textbf{Approach~2} respectively. In the final 
Section~\ref{last} we complete the picture by observing that these 
quantum Lie algebras are isomorphic to the 
one obtained from the bicovariant calculus 
on the Jordanian quantum group $SL_{h}(2)$.

\section{The Jordanian quantum groups}
The 2-parameter Jordanian quantum group $GL_{h,g}(2)$ is the 
co-quasitriangular Hopf algebra derived from the following $R$-matrix,
\begin{equation}
R=\begin{pmatrix}
1&-h&h&gh\\
0&1&0&-g\\
0&0&1&g\\
0&0&0&1\label{thermat}
\end{pmatrix}.
\end{equation}
Being triangular, i.e.\ $R_{21}R=I$, this $R$-matrix is trivially 
Hecke, with
\begin{equation}
(\hat{R}-1)(\hat{R}+1)=0,
\end{equation}
where $\hat{R}=PR$ and $P_{ij,kl}=\d_{il}\d_{jk}$. 

The quantum group associated with the $R$-matrix $R_{21}=R(-h,-g)$, with 
$g=h=1$, was first investigated by
Demidov et al.~\cite{DMMZ}, while $R_{21}$ with $g=h$ is the 1-parameter 
non-standard $R$-matrix whose quantum group was considered by
Zakrzewski~\cite{Zak}. Lazarev and Movshev\cite{LM} considered the 
quantum group associated with $R$ with $g=h$ and also the 
corresponding quantised universal enveloping algebra, $U_{h}(\slt)$. 
The quantised enveloping algebra was also investigated by 
Ohn~\cite{Ohn} and will be discussed in more detail in 
Section~\ref{defuhslt}.  In fact, the 
$R$-matrix,~(\ref{thermat}), can be extracted from an early work of 
Gurevich~\cite{Gur};  though the associated quantum group structure was 
not investigated there. 

In the usual way, defining an algebra valued matrix $T$ as
\begin{equation}
T=\begin{pmatrix}
a & b\\
c & d
\end{pmatrix},
\end{equation}
the relations of a matrix element bialgebra $A(R)$ are obtained from the 
well-known FRT~\cite{FRT} matrix relation,
\begin{equation}
RT_{1}T_{2}=T_{2}T_{1}R,
\end{equation}
as
\begin{gather}
ca=ac-gc^2,\;\;\;\;cd=dc-hc^2,\nonumber\\
db=bd+g(ad-bc+hac-d^2),\nonumber\\
ab=ba+h(ad-bc+hac-a^2),\nonumber\\
cb=bc-hac-gdc+ghc^2,\;\;\;\;da=ad+hac-gdc.\label{rels}
\end{gather}
The coalgebra structure is provided by a coproduct defined on the 
generators as,
\begin{gather}
\D(a)=a\ot a+b\ot c,\;\;\;\;\D(b)=a\ot b+b\ot d,\nonumber\\
\D(c)=c\ot a+d\ot c,\;\;\;\;\D(d)=c\ot b+d\ot d,\label{cps}
\end{gather}
with the counit given by,
\begin{gather}
\ep(a)=1,\;\;\;\;\ep(b)=0,\nonumber\\
\ep(c)=0,\;\;\;\;\ep(d)=1.\label{cous}
\end{gather}

$\hat{R}$ has a spectral decomposition,
\begin{equation}
\hat{R}=P^{+}+P^{-},
\end{equation}
where $P^{+}=\frac{1}{2}(\hat{R}+I)$ is a rank 3 projector and 
$P^{-}=\frac{1}{2}(\hat{R}-I)$ is a rank 1 projector. In the 
notation of Majid~\cite{Majbook}, these 
projectors provide the associated 
quantum co-plane, $\AA^{2|0}_{-1}$, and plane, $\AA^{2|0}_{1}$, respectively 
through the relations 
$P^{\pm}\mathbf{x}_{1}\mathbf{x}_{2}=0$ where $\mathbf{x}$ is the 
$2\times 1$ column vector 
$\bigl(\begin{smallmatrix}x_{1}\\x_{2}\end{smallmatrix}\bigr)$. 
These are associative algebras generated by the elements 
$x_{1}$ and $x_{2}$ subject to the relations,
\begin{equation}
x_{2}^{2}=0,\;\;\;\;x_{1}^{2}=-x_{2}x_{1},\;\;\;\;x_{1}x_{2}=-x_{2}x_{1},
\end{equation}
in the case of the co-plane $\AA^{2|0}_{-1}$, and,
\begin{equation}
x_{1}x_{2}=x_{2}x_{1}+gx_{2}^{2} \label{plane}
\end{equation}
in the case of the plane $\AA^{2|0}_{1}$. A result of 
Mukhin~\cite[Theroems~1,2,3]{Muk} (see also~\cite[Theorem~3.5]{Dem}) then tells us that 
$A(R)$ is the \emph{universal coacting 
bialgebra} on this pair of algebras in the sense of 
Manin~\cite{Man1}. In the language of 
Sudbery~\cite{Sud1}, $\AA^{2|0}_{-1}$ and $\AA^{2|0}_{1}$ are then \emph{complementary coordinate 
algebras} determining $A(R)$, and with an easy application of the Diamond 
Lemma~\cite{Berg} telling us that they are moreover \emph{superpolynomial 
algebras} generated by odd and even generators respectively having 
ordering algorithms with respect to the ordering of the generators,  
$x_{2}\prec x_{1}$, we deduce immediately from a result of 
Sudbery~\cite[Theorem~3]{Sud1} that $A(R)$ 
has as a basis the ordered monomials 
$\{b^{\a}a^{\b}d^{\g}c^{\d}:\a,\b,\g,\d\in\ZZ_{\geq 0}\}$. This fact 
is used extensively in the computations which lead to our main results.

It follows from $\AA^{2|0}_{-1}$, that $R$ is 
\emph{Frobenius}~\cite{C&Pbook} and so 
in the usual way we can obtain a group-like element in the bialgebra 
$A(R)$, $\mathcal{D}$, called the quantum determinant, and given by,
\begin{equation}
\mathcal{D}=ad-bc+hac.
\end{equation}
The commutation relations between $\mathcal{D}$ and the generators of $A(R)$ are,
\begin{gather}
\mathcal{D}a=a\mathcal{D}+(h-g)c\mathcal{D},\;\;\;\;
\mathcal{D}d=d\mathcal{D}-(h-g)c\mathcal{D},\nonumber \\
\mathcal{D}c=c\mathcal{D},\;\;\;\;
\mathcal{D}b=b\mathcal{D}+(h-g)(d\mathcal{D}-a\mathcal{D}-(h-g)c\mathcal{D}),
\end{gather}
so we can localise with respect to the Ore set 
$S=\{\mathcal{D}^{\a}:\a\in\ZZ_{\geq 1}\}$, and define  
$GL_{h,g}(2)=A(R)[\mathcal{D}^{-1}]$, having extra commutation relations,
\begin{gather}
a\mathcal{D}^{-1}=\mathcal{D}^{-1}a+(h-g)\mathcal{D}^{-1}c,\;\;\;\;
d\mathcal{D}^{-1}=\mathcal{D}^{-1}d-(h-g)\mathcal{D}^{-1}c,\nonumber\\
c\mathcal{D}^{-1}=\mathcal{D}^{-1}c,\;\;\;\;b\mathcal{D}^{-1}=\mathcal{D}^{-1}b+(h-g)(\mathcal{D}^{-1}d-\mathcal{D}^{-1}a-(h-g)\mathcal{D}^{-1}c),
\end{gather}
with
\begin{equation}
\D(\mathcal{D}^{-1})=\mathcal{D}^{-1}\ot \mathcal{D}^{-1},\;\;\;\;\ep(\mathcal{D}^{-1})=1.
\end{equation}
$GL_{h,g}(2)$ is a Hopf algebra with the antipode given by
\begin{gather}
S(a)=\mathcal{D}^{-1}(d+gc),\;\;\;\;S(b)=\mathcal{D}^{-1}(gd-ga-b+g^{2}c),\nonumber\\
S(c)=-\mathcal{D}^{-1}c,\;\;\;\;S(d)=\mathcal{D}^{-1}(a-gc),\;\;\;\;S(\mathcal{D}^{-1})=\mathcal{D}.
\end{gather}
The Hopf algebra $GL_{h,g}(2)$ is clearly still polynomial with basis 
$\{\mathcal{D}^{-\a}b^{\b}a^{\g}d^{\d}c^{\zeta}:\a,\b,\g,\d,\zeta\in\ZZ_{\geq 0}\}$.

With $g=h$, $\mathcal{D}$ is central and we can consistently set 
$\mathcal{D}=1$ and pass to the quantum group $SL_{h}(2)$. The 
relations for $SL_{h}(2)$ are just~(\ref{rels}), but with the 
combination $ad$ replaced wherever it appears by $bc-hac+1$ and also
the further relation $ad=bc-hac+1$. With $g=h=0$ we recover the 
classical group coordinate rings.
\label{Jorqugro}

\section{Review of Woronowicz's bicovariant differential calculus}

We begin with the basic definitions.

\bde
A \emph{first order differential calculus} over an algebra $A$ is 
a pair $(\G,\dif)$ such that:
\begin{enumerate}
\item
$\G$ is an $A$-bimodule, i.e.\ 
\begin{equation}
(a\w)b=a(\w b)
\end{equation}
for all $a,b\in A$, $\w\in \G$, where the left and right actions 
which make $\G$ respectively a left $A$-module and right $A$-module are written multiplicatively.
\item
$\dif$ is a linear map, $\dif:A\mapto \G$.
\item
For any $a,b\in A$, the Leibniz rule is satisfied, i.e.\ 
\begin{equation}
\dif(ab)=\dif(a)b+a\dif(b).
\end{equation}
\item
The bimodule $\G$, or `space of 1-forms', is spanned by elements of 
the form $a\dif b$, $a,b\in A$.
\end{enumerate}
\ede

\bre
Given two first order differential calculi over an algebra 
$A$, $(\G,\dif)$ and $(\G',\dif')$, we say that they are \emph{isomorphic} if 
there is a bimodule isomorphism $\lgh:\G\mapto\G'$ such that 
$\lgh\circ\dif=\dif'$. 
\label{iso}
\ere

\bre
We will usually write $\dif a$ for $\dif(a)$.
\ere

\bde
A \emph{bicovariant bimodule} over a Hopf algebra $A$ is a triple 
$(\G,\lca,\rca)$ such that:
\begin{enumerate}
\item
$\G$ is an $A$-bimodule. 
\item
$\G$ is an $A$-bicomodule with left and right coactions $\lca$ and 
$\rca$ respectively, i.e.\ 
\begin{equation}
(\id\ot\lca)\circ\lca=(\D\ot\id)\circ\lca,\;\;\;\;(\ep\ot\id)\circ\lca=\id, 
\label{lcact}
\end{equation}
making $\G$ a left $A$-comodule,
\begin{equation}
(\rca\ot\id)\circ\rca=(\id\ot\D)\circ\rca,\;\;\;\;(\id\ot\ep)\circ\rca=\id,
\label{rcact}
\end{equation}
making $\G$ a right $A$-comodule, and 
\begin{equation}
(\id\ot\rca)\circ\lca=(\lca\ot\id)\circ\rca,\label{bcom}
\end{equation}
which is the $A$-bicomodule property.
\item
The coactions, $\lca$ and $\rca$ are bimodule maps, i.e.\ 
\begin{gather}
\lca(a\w b)=\D(a)\lca(\w)\D(b),\label{lcabm}\\
\rca(a\w b)=\D(a)\rca(\w)\D(b).\label{rcabm}
\end{gather}
\end{enumerate}
\ede

\bre
The Sweedler notation for coproducts in the Hopf algebra $A$ is taken 
to be $\D(a)=a_{(1)}\ot a_{(2)}$ for all $a\in A$ and is extended to 
the coactions as $\lca(\w)=\w_{(A)}\ot \w_{(\G)}$ and 
$\rca(\w)=\w_{(\G)}\ot \w_{(A)}$. In this notation the 
conditions~(\ref{lcact}),~(\ref{rcact}) and~(\ref{bcom}) become,
\begin{gather}
\w_{(A)}\ot (\w_{(\G)})_{(A)}\ot (\w_{(\G)})_{(\G)}=(\w_{(A)})_{(1)}\ot 
(\w_{(A)})_{(2)}\ot \w_{(\G)},\;\;\;\; \ep(\w_{(A)})\w_{(\G)}=\w, \\
(\w_{(\G)})_{(\G)}\ot (\w_{(\G)})_{(A)}\ot \w_{(A)}=\w_{(\G)}\ot 
(\w_{(A)})_{(1)}\ot (\w_{(A)})_{(2)},\;\;\;\; 
\ep(\w_{(A)})\w_{(\G)}=\w, \\
\w_{(A)}\ot (\w_{(\G)})_{(\G)}\ot 
(\w_{(\G)})_{(A)}=(\w_{(\G)})_{(A)}\ot (\w_{(\G)})_{(\G)}\ot 
\w_{(A)},
\end{gather}
for all $\w\in\G$.
\ere 

\bde
A \emph{first order bicovariant differential calculus} over a Hopf 
algebra $A$ is a quadruple $(\G,\dif,\lca,\rca)$ such that:
\begin{enumerate}
\item
$(\G,\dif)$ is a first order differential calculus over $A$.
\item
$(\G,\lca,\rca)$ is a bicovariant bimodule over $A$.
\item
$\dif$ is both a left and a right comodule map, i.e.\ 
\begin{gather}
(\id\ot\dif)\circ\D(a)=\lca(\dif a),\label{dlcm}\\
(\dif\ot\id)\circ\D(a)=\rca(\dif a),\label{drcm}
\end{gather}
for all $a\in A$.
\end{enumerate}
\ede

\bre
Given a first order calculus over a Hopf algebra, $(\G,\dif)$,~(\ref{dlcm}) 
and~(\ref{drcm}) uniquely determine left and right coactions, and hence 
a bicovariant bimodule structure. However the \emph{existence} of 
these coactions and the corresponding bicovariant bimodule is of course 
\emph{not} guaranteed. 
\ere

\bex
Given any associative algebra $A$, we may form a first 
order differential calculus over $A$, $(A^{2},\udif)$, where
\begin{equation}
A^{2}=\left\{\sum_{k} a_{k}\ot b_{k}\in A\ot A:\sum_{k} a_{k}b_{k}=0\right\},
\end{equation}
and $\udif:A\mapto A^{2}$ is given by
\begin{equation}
\udif a=1\ot a-a\ot 1.
\end{equation}
$A^{2}$ has a bimodule structure given, for all $a_{k},b_{k},c\in A$ by 
\begin{equation}
c\left(\sum_{k} a_{k}\ot b_{k}\right)=\sum_{k} ca_{k}\ot b_{k},\;\;\;\;
\left(\sum_{k} a_{k}\ot b_{k}\right)c=\sum_{k} a_{k}\ot b_{k}c.
\end{equation}
The importance of this differential calculus lies in the fact that 
\emph{any} first order differential calculus over $A$, $(\G,\dif)$, is 
isomorphic to one of the form $(A^{2}/\mathcal{N},\pi\circ\udif)$ where 
$\mathcal{N}$ is the kernel of the surjective map $\pi:A^{2}\mapto\G$ 
defined by $\pi(\sum_{k} a_{k}\ot b_{k})=\sum a_{k}\dif b_{k}$. For this 
reason, $(A^{2},\udif)$ is said to be \emph{universal}. 
Moreover, it is not difficult to check that when $A$ is a Hopf algebra,
$(A^{2},\udif)$ is a first order bicovariant differential calculus.
\eex

\bde
An element $\w$ of a bicovariant bimodule $(\G,\lca,\rca)$ is said to 
be \emph{left-invariant} if
\begin{equation}
\lca(\w)=1\ot\w,\label{lica}
\end{equation}
\emph{right-invariant} if
\begin{equation}
\rca(\w)=\w\ot 1,
\end{equation}
and \emph{bi-invariant} if it is both left- and right-invariant.
\ede

\bre
Denoting the vector space of all left-invariant elements of $\G$ by 
$\G_{\inv}$, there is a projection $P:\G\mapto\G_{\inv}$ defined on 
any element $\w\in\G$ as,
\begin{equation}
P(\w)=S(\w_{(A)})\w_{(\G)}.
\end{equation}
If $(\G,\dif,\lca,\rca)$ is a first order bicovariant differential calculus 
over $A$ then it is not difficult to see the equivalence of the 
statements that the differentials generate $\G$ as a left $A$-module and 
that the elements $P(\dif a)=S(a_{(1)})\dif a_{(2)}$, for all $a\in A$, 
span the vector space of left-invariant 1-forms. Further, for any 
$a\in A$, $\dif a$ may be expressed in terms of left-invariant forms as
\begin{equation}
\dif a=a_{(1)}P(\dif a_{(2)}). \label{difininv}
\end{equation}
\label{lispan}
\ere

An alternative construction of the universal differential calculus 
occurs in the situation of particular interest to us. It is described 
in the following example.

\bex
We start with a Hopf algebra 
$A$ defined in terms of a finite number of generators, 
$(T_{ij})_{i,j=1\ldots n}$ say, 
together with certain relations which are consistent with a PBW-type basis and the 
usual (matrix-element bialgebra) coproduct and counit. We may 
then introduce a bimodule $\G_{0}$ as the free $A$-bimodule on the 
symbols $\dif_{0} T_{ij}$ with the linear map $\dif_{0}:A\mapto\G$ defined on 
any element of $A$ by way of the Leibniz rule. $(\G_{0},\dif_{0})$ is then 
also a universal first order differential calculus, and being 
therefore
isomorphic to $(A^{2},\udif)$ is certainly bicovariant. 
Indeed,~(\ref{dlcm}) and~(\ref{drcm}) specify the coactions, and the space of 
left invariant forms, ${\G_{0}}_{\inv}$, is spanned by all elements of 
the form $S(a_{(1)})\t_{ik}a_{(2)}$ where $\t_{ik}=\sum_{j} 
S(T_{ij})\dif_{0}T_{jk}$ and $a\in A$. If we denote by 
$(\vartheta_{i})_{i\in I}$ a basis for $\G_{0}$, where $I$ is  some 
countably infinite index set, then from~(\ref{drcm}), there exist 
$v_{ij}$s such that the right coaction on the 
$\vartheta_{i}$s takes the form,
\begin{equation}
\rca(\vartheta_{i})=\sum_{j\in I}\vartheta_{j}\ot v_{ji}.
\end{equation}
In particular, the $n^{2}$ left invariant elements $\t_{ik}$ span a 
sub-bicomodule, with
\begin{equation}
\rca(\t_{ik})=\sum_{s,t=1\ldots n}\t_{st}\ot 
S(T_{is})T_{tk}.\label{univv}
\end{equation}
\label{relexamp}
\eex

At this point in Example~\ref{relexamp} 
the only relations between algebra and bimodule elements 
are those coming from the Leibniz rule. Obtaining further relations 
involves finding a suitable 
`relation space', $\mathcal{N}$, which can be factored out from 
$\G_{0}$ while maintaining the bicovariance. In the context of the 
universal calculus, $(A^{2},\udif)$, Theorems~1.5 and~1.8 
of~\cite{Wor} tell us that such $\mathcal{N}$ must be of the form 
$\tau^{-1}(A\ot\mathcal{R})$ where $\tau^{-1}(a\ot b)=aS(b_{(1)})\ot 
b_{(2)}$ and $\mathcal{R}$ is a right ideal of $A$, contained in $\ker\ep$, 
which is stable under 
the right adjoint coaction\footnote{We recall that the right adjoint 
coaction is defined on any $a\in A$ as $Ad^{*}_{R}(a)=a_{(2)}\ot 
S(a_{(1)})a_{(3)}$}. Conversely, given a first order bicovariant 
differential calculus, $(\G,\dif,\lca,\rca)$, $\mathcal{R}$ may be 
recovered as the set of all $a\in\ker\ep$ such that $P(\dif a)=0$.

It is desirable to classify all bicovariant calculi on a given 
quantum group which have particular properties. This problem of 
classification can be regarded as the problem 
of classifying the ad-invariant ideals $\mathcal{R}$~\cite{SS2,SS3,SS4}.
However, in this paper, 
following M\"{u}ller-Hoissen~\cite{MH1}, we 
will consider a more `hands on' approach. We look for calculi 
whose 
bimodule of 1-forms is generated as a \emph{left} $A$-module by the 
differentials of the generators (this assumption is also made 
in~\cite{SS2,SS3}). In effect we are passing directly to 
a class of, as yet unspecified, quotients of $\G_{0}$ which we then 
wish to constrain with by the requirement that the bicovariance is not 
destroyed. For 
this approach we need the characterisation of bicovariant bimodules 
which is provided by the following theorem of Woronowicz.

\bth
Let $(\G,\lca,\rca)$ be a bicovariant bimodule over $A$ and let 
$(\t_{i})_{i\in I}$ be a basis of $\G_{\inv}$, where $I$ is some 
countable index set. Then we have:
\begin{enumerate}
\item
Any element $\w\in\G$ has a unique expression as
\begin{equation}
\w=\sum_{i\in I} a_{i}\t_{i},
\end{equation}
where $a_{i}\in A$.
\item
There exist linear functionals $f_{ij}\in A^{*}$, $i,j\in I$, such that
\begin{equation}
\t_{i}a=\sum_{j\in I} (f_{ij}\star a)\t_{j}, \label{thfcrs}
\end{equation}
where $\a\star a=\a(a_{(2)})a_{(1)}$ for all $\a\in A^{*}$, $a\in A$.
\item
The functionals $f_{ij}$ are uniquely determined by~(\ref{thfcrs}) and 
satisfy
\begin{equation}
f_{ik}(ab)=\sum_{j\in I} f_{ij}(a)f_{jk}(b),\;\;\;\;f_{ik}(1)=\d_{ik}, \label{frep}
\end{equation}
for all $a,b\in A$, $i,j\in I$.
\item
There exist elements $v_{ij}\in A$, $i,j\in I$, such that for all 
$t_{i}$,
\begin{gather}
\rca(\t_{i})=\sum_{j\in I} \t_{j}\ot v_{ji},\label{adrep}\\
\D(v_{ik})=\sum_{j\in I} v_{ij}\ot 
v_{jk},\;\;\;\;\ep(v_{ik})=\d_{ik}.\label{rrep}
\end{gather}
\item
For all $a\in A$,
\begin{equation}
\sum_{j\in I} v_{ji}(a\star f_{jk})=\sum_{j\in I} (f_{ij}\star a)v_{kj},\label{lrtt}
\end{equation}
where $a\star\a=\a(a_{(1)})a_{(2)}$ for all $\a\in A^{*}$, $a\in A$.
\end{enumerate}
Conversely, if $(\t_{i})_{i\in I}$ is a basis of a vector space $V$, 
and we have functionals $(f_{ij})_{i,j\in I}$ defined on $A$ and 
elements $(v_{ij})_{i,j\in I}$ in $A$ which satisfy~(\ref{frep}), 
(\ref{rrep}) and~(\ref{lrtt}), then there exists a unique bicovariant 
bimodule such that $V=\G_{inv}$ and~(\ref{thfcrs}) and~(\ref{adrep}) 
are satisfied.
\label{basthm}
\ethe

This result places significant constraints on the possible bicovariant 
calculi which are consistent with our assumption that the 
differentials of the generators should generate the bimodule of forms 
as a left $A$-module. That 
assumption implies immediately that $\G_{\inv}$ is spanned by 
the \emph{finite} set of elements $\t_{ik}$. Choosing a basis from 
this set we either take all $n^{2}$ as linear independent elements, in 
which case the $v_{ij}$s of~(\ref{adrep}) and~(\ref{rrep}) have 
already been determined in~(\ref{univv}), or we introduce linear relations between 
the $\t_{ik}$s. In the latter case we must be sure that such 
relations don't destroy the bicovariance --- they must factor 
through the coactions. This condition will tend to fix the possible 
bases of $\G_{\inv}$, which in turn determines the $v_{ij}$s 
through~(\ref{adrep}) and~(\ref{univv}). It then follows immediately that~(\ref{rrep}) 
is satisfied. Now, trying to introduce commutation relations which factor through 
the left and right coactions which we wish to maintain, it is not too 
difficult to see that when these relations take the 
form~(\ref{thfcrs}) with the $f_{ij}$s satisfying~(\ref{frep}) 
and~(\ref{lrtt}), bicovariance is maintained. Moreover, as the theorem 
states, the commutation relations in \emph{any} bicovariant bimodule 
must \emph{always} take this form. Thus our method of attacking the 
classification problem can now be discerned. We decide upon a valid basis of 
$\G_{inv}$, and then assume the general form for the commutation 
relations,~(\ref{thfcrs}). We must then impose as constraints, consistency with the relations 
already present from the Leibniz rule, together with~(\ref{frep}) and~(\ref{lrtt}). 
This procedure will be made absolutely explicit for the case of the 
quantum groups of particular interest to us in the following section.

\bre
Having chosen a basis $(\t_{i})_{i\in I}$ for $\G_{\inv}$, the 
right-invariant
elements $(\eta_{i})_{i\in I}$ defined by $\eta_{i}=\sum_{j\in I} 
\t_{j}S(v_{ji})$ form a basis for the vector space of right invariant 
1-forms, $\G_{\rinv}$.
\ere

\bre
The \emph{dimension} of a bicovariant calculus is defined to be $\dim 
\G_{\inv}$. We will only be interested in finite, $d$-dimensional, 
examples so we will eschew the index set and consider indices 
running over a finite set.
\ere

Given two bicovariant bimodules, $(\G,\lca,\rca)$ and 
$(\tilde{\G},\lcad,\rcad)$, of dimensions $d$ and $\tilde{d}$ 
respectively, their tensor product over $A$, 
$(\G'=\G\ot_{A}\tilde{\G},\lcat,\rcat)$ is also a 
bicovariant bimodule as follows. The left and right coactions of $A$ 
on $\G'$ are given by
\begin{equation}
\lcat(\w\ot\tilde{\w})={\w}_{(A)}{\tilde{\w}}_{(A)}\ot{\w}_{(\G)}\ot{\tilde{\w}}_{(\G)},\;\;\;\;
\rcat(\w\ot\tilde{\w})={\w}_{(\G)}\ot{\tilde{\w}}_{(\G)}\ot{\w}_{(A)}{\tilde{\w}}_{(A)},
\end{equation}
and the $f$ functionals of Theorem~\ref{basthm} are now given by 
$F'_{ik,jl}={f}_{ij}*{\tilde{f}}_{kl}$, where $*$ is the usual convolution 
product, such that
\begin{equation}
({\t}_{i}\ot{\tilde{\t}}_{k})a=\sum_{j=1\ldots d, l=1\ldots\tilde{d}}(F'_{ik,jl}\star 
a)({\t}_{j}\ot{\tilde{\t}}_{l}),
\end{equation}
for all $a\in A$, where $({\t}_{i})_{i=1\ldots d}$ and 
$({\tilde{\t}}_{i})_{i=1\ldots 
\tilde{d}}$ are the bases of left invariant elements in ${\G}_{\inv}$ and 
${\tilde{\G}}_{\inv}$ respectively, so that 
$({\t}_{i}\ot{\tilde{\t}}_{j})_{i=1\ldots d, j=1\ldots\tilde{d}}$ is the basis of left 
invariant elements in ${\G'}_{\inv}$.
It is clear that in this way we can build arbitrary n-fold tensor powers, 
$\G^{\ot n}=\G\ot_{A}\G\ot_{A}\ldots\ot_{A}\G$ 
of a given bicovariant bimodule, all of which are themselves 
bicovariant with left and right coactions denoted $\lcan$ and $\rcan$ 
respectively. We can then define $\G^{\ot}=A\opl\G\opl\G^{\ot 
2}\opl\ldots$ to be the analog of the classical algebra of covariant 
tensor fields. $\G^{\ot}$ is a bicovariant graded algebra, 
in that it is a tensor algebra and a bicovariant bimodule over $A$, with  
coactions $\lcata$ and $\rcata$ which are algebra maps and coincide 
on elements of $\G^{\ot n}$ with the coactions $\lcan$ and $\rcan$ 
respectively.

The next step in Woronowicz's construction of a noncommutative 
geometry for Hopf algebras is to introduce an analogue of the 
classical external algebra of forms. Starting from the bicovariant 
graded algebra $\G^{\ot}$ we want to obtain another 
bicovariant graded algebra, the \emph{external bicovariant graded 
algebra}, $\EQ$, as a quotient, $\EQ=\G^{\ot}/S$, by some 
graded two-sided ideal $S$. Woronowicz introduces the \emph{unique} linear bimodule 
map $\L:\G\ot\G\mapto\G\ot\G$ such that 
\begin{equation}
\L(\t\ot\eta)=\eta\ot\t \label{defbraid}
\end{equation} 
for any $\t\in\G_{\inv}$, $\eta\in\G_{\rinv}$. $\L$ is then also a 
bicomodule map, is given 
explicitly on the basis of left invariant elements 
$(\t_{i}\ot\t_{k})_{i,j=1\ldots d}$ by
\begin{equation}
\L(\t_{i}\ot\t_{k})=\sum_{s,t=1\ldots 
d}\L_{ik,st}\t_{s}\ot\t_{t}=\sum_{s,t=1\ldots 
d}f_{it}(v_{sk})\t_{s}\ot\t_{t},\label{braidcomp}
\end{equation}
and can be shown to satisfy the braid equation, 
$\L_{12}\L_{23}\L_{12}=\L_{23}\L_{12}\L_{23}$. The map $\L$ may then be 
used in just the same way as the permutation operator is used 
classically. Thus we define an analogue of the antisymmetrisation operator 
on $\G^{\ot n}$, $W_{1\ldots n}$, by replacing the classical permutation 
operator by $\L$ everywhere in the classical antisymmetriser. $W_{1\ldots 
n}$ is then a bimodule and bicomodule map $W_{1\ldots n}:\G^{\ot 
n}\mapto\G^{\ot n}$ and we can define $S^{n}=\ker W_{1\ldots n}$ so 
that $\G^{\ot n}/S^{n}$ is isomorphic to $\Img W_{1\ldots n}$ and 
$\EQ=\G^{\ot}/S$ where $S=\bigoplus_{n=2} S^{n}$. As $W_{1\ldots n}$ 
is a bicomodule map, the left and right coactions of $\G^{\ot}$ descend 
to the quotient where they shall be denoted $\lcae$ and $\rcae$ 
respectively. Moreover, we can now 
define the \emph{wedge product} as 
$\w_{1}\wedge\w_{2}\wedge\ldots\wedge\w_{n}=W_{1\ldots 
n}(\w_{1}\ot\w_{2}\ot\ldots\ot\w_{n})$ where the $\w_{i}$ are arbitrary 
elements of $\G$.

\bre
If we were to construct the bicovariant graded algebra 
${\tilde{\G}}^{\ot}$ from a bicovariant bimodule ${\tilde{\G}}$ which 
contains $\G$ as a sub-bimodule and a sub-bicomodule , then the 
bicovariant graded algebra $\G^{\ot}$ embeds naturally into 
${\tilde{\G}}^{\ot}$ in the sense that there is the obvious natural embedding 
map $\phi:\G^{\ot}\mapto{\tilde{\G}}^{\ot}$ such that
\begin{gather}
\phi|_{A}=\id,\;\;\;\;\phi|_{\G}=\iota,\\
\phi(\w\rho)=\phi(\w)\phi(\rho),\\
(\id\ot\phi)\circ\lca=\lcata\circ\phi,\;\;\;\;(\phi\ot\id)\circ\rca=\rcata\circ\phi,
\end{gather}
where $\iota:\G\mapto\tilde{\G}$ is the natural inclusion. An attractive feature of 
the external bicovariant algebra 
construction, which has its origin in the uniqueness of the map $\L$, is that 
this embedding survives the quotienting 
procedures so that $(\EQ,\lcae,\rcae)$ embeds naturally in 
$(\tilde{\EQ},\lcaed,\rcaed)$.
\ere

We may now state the following theorem of Woronowicz.

\bth
Let $\EQ$ be the external bicovariant algebra constructed above. There 
exists one and only one linear map $\dif:\EQ\mapto\EQ$ such that:
\begin{enumerate}
\item
$\dif$ increases the grade by one.
\item
On elements of grade 0, $\dif$ coincides with the differential of the 
first order bicovariant calculus $(\G,\dif,\lca,\rca)$.
\item
For all $\w\in\G^{\ot k}$, $k=0,1,2,\ldots$ and $\w'\in\EQ$,
\begin{equation}
\dif(\w\wedge\w')=\dif\w\wedge\w'+(-1)^{k}\w\wedge\dif\w'.\label{hodif}
\end{equation}
\item
For any $\w\in\EQ$,
\begin{equation}
\dif(\dif\w)=0.\label{dsqr}
\end{equation}
\item
$\dif$ is a both a left and a right comodule map i.e.\
\begin{gather}
\lcae(\dif\w)=(\id\ot\dif)\circ\lcae(\w),\\
\rcae(\dif\w)=(\dif\ot\id)\circ\rcae(\w),
\end{gather}
for all $\w\in\EQ$.
\end{enumerate}
\label{hothm}
\ethe

\bre
The external bicovariant algebra $\EQ$ equipped with the differential 
described in this result will be called the \emph{exterior bicovariant 
differential calculus} over the Hopf algebra $A$, and denoted 
$(\EQ,\dif,\lcae,\rcae)$. Brzezi\'{n}ski has shown that this is a 
super-Hopf algebra~\cite{Brz}.
\ere

\bre
To prove this theorem Woronowicz extends the bimodule $\G$ to 
$\tilde{\G}=Ax\opl\G$ where $Ax$ is the left $A$-module freely 
generated by the single element $x$, such that the right action of 
$A$ on an element $ax$ is given by $axb=abx+a\dif b$ and the 
coactions are such that the element $x$ is bi-invariant. The theorem is then 
proved for the external bicovariant algebra $\tilde{\EQ}$ built from $\tilde{\G}$ 
with the final result for $\EQ$ coming after using the natural 
embedding of $\EQ$ in $\tilde{\EQ}$. Along the way the 
differential is expressed as $\dif a=[x,a]$, for any $a\in A$, and more generally 
as $\dif \w=[x,\w]_{\mp}$ for any $\w\in\tilde{\G}$ where 
$[x,\w]_{\mp}=x\wedge\w\mp\w\wedge x$ with $-$ and $+$ for $\w$ of 
even and odd grade respectively. In particular examples of bicovariant 
differential calculi it sometimes happens that there is a bi-invariant 
element \emph{within} the unextended calculus which implements the 
differential in this way. Such calculi are called \emph{inner}. 
\label{reminner}
\ere

The final component of the Woronowicz differential calculus is the 
analogue of the classical Lie algebra of tangent vectors at the 
identity. Classically this is isomorphic to 
$(\ker\ep/(\ker\ep)^{2})^{*}$. In the abstract Hopf algebra setting it 
turns out that defining a space $\mathcal{L}$ by $\mathcal{L}=(\ker\ep/\mathcal{R})^{*}$ 
there is a unique bilinear map, $\pair{}{}:\G\times\mathcal{L}\mapto\CC$, such 
that for all $a\in A$, $\w\in\G$ and $\chi\in\mathcal{L}$,
\begin{equation}
\pair{a\w}{\chi}=\ep(a)\pair{\w}{\chi},\;\;\;\;\pair{\dif a}{\chi}=\chi(a),
\end{equation}
which is \emph{non-degenerate} as a pairing between $\G_{\inv}$ and $T$. 
Thus, we can 
introduce a basis $(\chi_{i})_{i=1\ldots d}$ for $\mathcal{L}$ dual to 
$(\t_{i})_{i=1\ldots d}$. Also, it is not too difficult to see that for any 
$\w\in\G$ there is an element $a\in\ker\ep$ such that $P(\w)=P(\dif 
a)$. Combined with the defining characteristics of the pairing 
between $\G$ and $\mathcal{L}$, we then arrive at a basis $a_{i}$ of 
$(\ker\ep/\mathcal{R})$ such that $\t_{i}=P(\dif a_{i})$. These 
$a_{i}$ are analogues of the classical coordinate functions at the 
identity. The following result of Woronowicz continues the analogy 
with the classical situation.

\bth
For all $a,b\in A$,
\begin{gather}
\dif a=\sum_{j=1\ldots d}(\chi_{j}\star a)\t_{j},\\
\dif \t_{i}=-\sum_{j,l=1\ldots 
d}\mathscr{C}_{jl,i}\t_{j}\wedge\t_{l},\label{cmeq}\\
\chi_{i}(ab)=\sum_{j=1\ldots d}\chi_{j}(a)f_{ji}(b)+\ep(a)\chi_{i}(b),\label{ptoc}\\
\sum_{j=1\ldots 
d}\chi_{j}(a)v_{ij}=\chi_{i}(a_{(2)})S(a_{(1)})a_{(3)},\label{ifb}
\end{gather}
where the $f_{ji}$ and $v_{ij}$ are as introduced in 
Theorem~\ref{basthm}, 
and $\mathscr{C}_{jl,i}=(\chi_{j}*\chi_{l})(a_{i})$.
\label{lathm}
\ethe

The equation~(\ref{ptoc}) shows us that the $\chi_{i}$ may be 
interpreted as `deformed derivations' and is equivalent to a coproduct for the 
`quantum tangent vectors',
\begin{equation}
\D(\chi_{i})=\sum_{j=1\ldots d}\chi_{j}\ot f_{ji}+1\ot\chi_{i}.
\end{equation}
Further, as $\mathcal{R}$ is stable under the right adjoint coaction, 
it follows that $\mathcal{L}$ is stable under the right adjoint 
$A^{*}$-action, $\radact$, that is, 
$S(\a_{(1)})\chi\a_{(2)}\in\mathcal{L}$ for all $\a\in A^{*}$ and 
any $\chi\in\mathcal{L}$. In the classical case this action, restricted to 
the tangent space, provides the Lie bracket. So we may define a 
\emph{quantum Lie bracket} as,
\begin{equation}
[\chi_{i},\chi_{k}]=\chi_{i}\radact\chi_{k}=\sum_{j=1\ldots 
d}\mathcal{C}_{ik,j}\chi_{j},\label{rhs1}
\end{equation}
where the $\mathcal{C}_{ik,j}$ are analogues of the classical Lie 
algebra structure constants and are still to 
be determined. But from the coproduct on the $\chi_{i}$ we can expand 
the right adjoint $A^{*}$-action, to obtain,
\begin{equation}
[\chi_{i},\chi_{k}]=\chi_{i}\chi_{k}-\sum_{s=1\ldots 
d}\chi_{s}(\chi_{i}\radact f_{sk}).
\end{equation}
We can then use~(\ref{ifb}) to determine $\chi_{i}\radact f_{sk}$ and 
obtain an expression for the bracket as a \emph{quantum commutator},
\begin{equation}
[\chi_{i},\chi_{k}]=\chi_{i}\chi_{k}-\sum_{s,t=1\ldots 
d}\L_{st,ik}\chi_{s}\chi_{t}.\label{rhs2}
\end{equation}
The structure constants, $\mathcal{C}_{ik,j}$ may now be determined by 
evaluating the right hand sides of~(\ref{rhs1}) and~(\ref{rhs2}) on 
$a_{l}$ and equating the results to obtain,
\begin{equation}
\mathcal{C}_{ik,j}=\mathscr{C}_{ik,j}-\sum_{s,t=1\ldots 
d}\L_{st,ik}\mathscr{C}_{st,j}.\label{qstrcon}
\end{equation}

There is a \emph{quantum Jacobi identity} for the quantum Lie bracket 
given by
\begin{equation}
[\chi_{i},[\chi_{j},\chi_{k}]]=[[\chi_{i},\chi_{j}],\chi_{k}]-
\sum_{s,t=1\ldots d}\L_{st,jk}[[\chi_{i},\chi_{s}],\chi_{t}].
\end{equation}

The \emph{universal enveloping algebra} of the quantum Lie algebra 
$\mathcal{L}$ may be introduced as the quotient of the tensor 
algebra of $\mathcal{L}$ by the 
two-sided ideal generated by the elements $\chi_{i}\chi_{k}-\sum_{s,t=1\ldots 
d}\L_{st,ik}\chi_{s}\chi_{t}-\sum_{j=1\ldots 
d}\mathcal{C}_{ik,j}\chi_{j}$.

\label{WorThr}

\section{The classification procedure}
In this section we make explicit the procedure, already outlined in the 
previous section, for determining under certain 
assumptions all possible first order bicovariant differential 
calculi for a given quantum group. It was first applied by M\"{u}ller-Hoissen~\cite{MH1,MH2} 
to the case of the standard 2-parameter quantum group $GL_{q,p}(2)$. Some further results 
appeared in subsequent papers~\cite{MH3,MHR}, and 
in~\cite{Bres} it was applied to the standard 1-parameter quantum 
group $GL_{q}(3)$. Here we apply this `recipe' to the cases of 
$GL_{h,g}(2)$ and $SL_{h}(2)$. 

Starting with our quantum group $A$, 
where $A$ is either $GL_{h,g}(2)$ or $SL_{h}(2)$, we 
introduce the 
first order differential calculus, in the first instance, as the free 
$A$-bimodule $\G_{0}$, on the 
symbols $\{\dif_{0} a,\dif_{0} b,\dif_{0} c,\dif_{0} d\}$ with the differential 
$\dif_{0}:GL_{h,g}(2)\mapto\G_{0}$ defined on 
any element of $A$ by way of the Leibniz rule. Note 
that $\dif_{0} \mathcal{D}^{-1}=\mathcal{D}^{-1}\dif \mathcal{D} \mathcal{D}^{-1}$. 
However, as already 
mentioned, we pass directly to some quotient $(\G,\dif)$ which we assume 
is generated as a left $A$-module by $\{\dif a,\dif b,\dif c,\dif 
d\}$ and is still a bicovariant bimodule, denoted $(\G,\dif,\lca,\rca)$. 
Then, by Remark~\ref{lispan} $\G_{\inv}$ is 
spanned by the 4 left-invariant forms 
$\{\t_{1},\t_{2},\t_{3},\t_{4}\}$ where,
\begin{gather}
\t_{1}=P(\dif a)=S(a)\dif a+S(b)\dif c,\;\;\;\;\t_{2}=P(\dif b)=S(a)\dif 
b+S(b)\dif d,\nonumber\\
\t_{3}=P(\dif c)=S(c)\dif a+S(d)\dif c,\;\;\;\;\t_{4}=P(\dif d)=S(c)\dif 
b+S(d)\dif d.\label{thtodif}
\end{gather}
In the other direction, from~(\ref{difininv}),
the differentials of the generators may be written in terms of the 
left invariant $\t_{i}$s as,
\begin{gather}
\dif a=a_{(1)}P(\dif a_{(2)})=a\t_{1}+b\t_{3},\;\;\;\;\dif 
b=b_{(1)}P(\dif b_{(2)})=a\t_{2}+b\t_{4},\nonumber\\
\dif c=c_{(1)}P(\dif c_{(2)})=c\t_{1}+d\t_{3},\;\;\;\;\dif 
d=d_{(1)}P(\dif d_{(2)})=c\t_{2}+d\t_{4}.\label{diftoth}
\end{gather}
We are now free to choose from, $\{\t_{1},\t_{2},\t_{3},\t_{4}\}$, 
a basis for $\G_{\inv}$ and look 
for corresponding non-trivial first order bicovariant differential 
calculi. We choose to look for calculi with the classical dimension, so for 
$GL_{h,g}(2)$ we will make the further assumption that 
$\{\t_{1},\t_{2},\t_{3},\t_{4}\}$ is a basis of $\G_{\inv}$ and call 
any extant calculi 4D calculi,
while for $SL_{h}(2)$ we will look for 3-dimensional, 3D, calculi by assuming 
that $\{\t_{1},\t_{2},\t_{3}\}$ 
is a basis of $G_{\inv}$ with 
$\t_{4}=\a_{1}\t_{1}+\a_{2}\t_{2}+\a_{3}\t_{3}$ and the coefficients 
$\a_{1}$, $\a_{2}$ and $\a_{3}$ to be determined. 
The assumptions will be justified if we
find non-trivial calculi. This in turn will be established if we can 
find functionals $f_{ij}$ and elements $v_{ij}$ as in Theorem~\ref{basthm} 
consistent with our assumptions. 

As we discussed in the previous 
section, the elements $v_{ij}$ are 
already fixed through our assumptions, by 
~(\ref{drcm}),~(\ref{thtodif}) and~(\ref{adrep}). 
For the 4D, $GL_{h,g}(2)$, calculi, we have
\begin{equation}
V=||v_{ij}||=\begin{pmatrix}
S(a)a & S(a)b & S(c)a & S(c)b \\
S(a)c & S(a)d & S(c)c & S(c)d \\
S(b)a & S(b)b & S(d)a & S(d)b \\
S(b)c & S(b)d & S(d)c & S(d)d
\end{pmatrix}.
\end{equation}
For the 3D, $SL_{h}(2)$, calculi we have
\begin{equation}
V=||v_{ij}||=\begin{pmatrix}
S(a)a+\a_{1}S(b)c & S(a)b+\a_{1}S(b)d & S(c)a+\a_{1}S(d)c \\
S(a)c+\a_{2}S(b)c & S(a)d+\a_{2}S(b)d & S(c)c+\a_{2}S(d)c \\
S(b)a+\a_{3}S(b)c & S(b)b+\a_{3}S(b)d & S(d)a+\a_{3}S(d)c
\end{pmatrix},
\end{equation}
but in this case we must also have,
\begin{equation}
\rca(\t_{4}-\a_{1}\t_{1}+\a_{2}\t_{2}+\a_{3}\t_{3})=0,\label{extsl2}
\end{equation}
Evaluating~(\ref{extsl2}) with~(\ref{thtodif}),~(\ref{drcm}) and~(\ref{diftoth}) 
fixes the $\a_{i}$ coefficients to be given by $\a_{1}=-1$, $\a_{2}=0$ 
and $\a_{3}=-2h$, so that we must have
\begin{equation}
\t_{4}=-\t_{1}-2h\t_{3}.\label{3Dred}
\end{equation}
In fact, when we look for bi-invariant forms in the 4D calculi, 
we soon find that up to scalar multiplication there is but one, 
which we will denote by $\Tr_{h}\Theta$, and which is given by
\begin{equation}
\Tr_{h}\Theta=\t_{1}+2h\t_{3}+\t_{4}
\end{equation}

Existence of the bicovariant first order differential calculi which we 
seek now hinges entirely on the $f_{ij}$s. Following 
M\"{u}ller-Hoissen~\cite{MH1}, let us refine our notation slightly. 
For the 4D calculi, we write the relations~(\ref{thfcrs}) in 
terms of four $4\times 4$ matrices, the `$ABCD$' matrices, $A_{ij}=f_{ij}(a)$, 
$B_{ij}=f_{ij}(b)$, $C_{ij}=f_{ij}(c)$ and $D_{ij}=f_{ij}(d)$, as
\begin{gather}
\t_{i}a=\sum_{j=1\ldots 
4}(aA_{ij}+bC_{ij})\t_{j},\;\;\;\;\t_{i}b=\sum_{j=1\ldots 
4}(aB_{ij}+bD_{ij})\t_{j},\nonumber\\
\t_{i}c=\sum_{j=1\ldots 
4}(cA_{ij}+dC_{ij})t_{j},\;\;\;\;\t_{i}d=\sum_{j=1\ldots 
4}(cB_{ij}+dD_{ij})\t_{j}.\label{thfmcrs}
\end{gather}
In the case of the 3D calculi the $4\times 4$ $ABCD$ matrices are simply 
replaced by $3\times 3$ $ABCD$ matrices with the summations then to 3.
We may now list the constraints on the $ABCD$ matrices.
\begin{description}
\item[Constraint~1]
Differentiating the quantum group relations,~(\ref{rels}), 
to obtain, in $R$-matrix form,
\begin{equation}
\dif(RT_{1}T_{2}-T_{2}T_{1}R)=R(\dif T_{1})T_{2}+RT_{1}(\dif 
T_{2})-(\dif T_{2})T_{1}R-T_{2}(\dif T_{1})R=0,
\end{equation}
we replace the differentials by left-invariant forms 
through~(\ref{diftoth}). We then 
use~(\ref{thfmcrs}) to commute the $\t_{i}s$ to 
the right which allows us then to equate the (ordered) algebra valued 
coefficients and obtain linear relations between the matrix elements of the 
$ABCD$ matrices. Similarly for the 3D case, but then, when replacing 
the differentials by left-invariant forms, we also use~(\ref{3Dred}).
\item[Constraint~2]
In both the 3D and 4D cases the relations~(\ref{lrtt}) may be expressed in 
the following matrix form,
\begin{equation}
\begin{pmatrix}
V^{T}A & V^{T}B\\
V^{T}C & V^{T}D
\end{pmatrix}
\begin{pmatrix}
a & b\\
c & d
\end{pmatrix}
=\begin{pmatrix}
a & b\\
c & d
\end{pmatrix}
\begin{pmatrix}
AV^{T} & BV^{T}\\
CV^{T} & DV^{T}
\end{pmatrix}.
\end{equation}
Only algebra elements appear in these equations and once all terms are 
`straightened' they yield further linear relations among the matrix elements of 
the $ABCD$ matrices.
\item[Constraint~3]
The equations~(\ref{frep}) are telling us that $A$, $B$, $C$ and 
$D$ must be the representation matrices of $a$, $b$, $c$ and $d$ 
respectively, and that the matrix representation of the determinant, 
$\DD$ say, where $\DD=AD-BC+hAC$, must be invertible in the case of $GL_{h,g}(2)$, and 
equal to the identity in the case of $SL_{h}(2)$. Imposing these 
conditions we obtain non-linear relations amongst the $ABCD$ matrix 
elements.
\end{description}

\bre
Recall from Remark~\ref{iso} our definition of isomorphism for differential 
calculi. As we assume a single $\t_{i}$ basis, different possible $ABCD$
matrices must correspond to non-isomorphic calculi.
\ere 

The $ABCD$ matrices which result from this procedure provide the most 
general possible first order bicovariant calculi under the stated 
assumptions. We may now investigate the external bicovariant graded 
algebras, and also the `quantum Lie algebras' which are related to 
our first order calculi. 

We begin by using Theorem~\ref{hothm}, in particular~(\ref{dsqr}), together 
with~(\ref{diftoth}) to deduce the 
structure constants $\mathscr{C}_{ij,k}$ appearing in the `Cartan-Maurer 
equations',~(\ref{cmeq}). For the 4D calculi, we obtain four $4\times 
4$ matrices,
\begin{align}
\mathscr{C}_{ij,1}&=\begin{pmatrix}
1&0&0&0\\
0&0&1&0\\
0&0&0&0\\
0&0&0&0
\end{pmatrix},&
\mathscr{C}_{ij,2}&=\begin{pmatrix}
0&1&0&0\\
0&0&0&1\\
0&0&0&0\\
0&0&0&0
\end{pmatrix},\nonumber\\
\mathscr{C}_{ij,3}&=\begin{pmatrix}
0&0&0&0\\
0&0&0&0\\
1&0&0&0\\
0&0&1&0
\end{pmatrix},&
\mathscr{C}_{ij,4}&=\begin{pmatrix}
0&0&0&0\\
0&0&0&0\\
0&1&0&0\\
0&0&0&1
\end{pmatrix},
\end{align}
while for the 3D calculi we obtain three $3\times 3$ matrices,
\begin{equation}
\mathscr{C}_{ij,3}=\begin{pmatrix}
1&0&0\\
0&0&1\\
0&0&0
\end{pmatrix},\;\;\;\;
\mathscr{C}_{ij,2}=\begin{pmatrix}
0&1&0\\
-1&0&-2h\\
0&0&0
\end{pmatrix},\;\;\;\;
\mathscr{C}_{ij,3}=\begin{pmatrix}
0&0&-1\\
0&0&0\\
1&0&-2h
\end{pmatrix},
\end{equation}
together with a relation,
\begin{equation}
2\t_{1}\wedge\t_{1}+4h\t_{3}\wedge\t_{1}+\t_{2}\wedge\t_{3}+\t_{3}\wedge\t_{2}=0,
\end{equation}
which will have to be consistent with any commutation relations which 
we derive for left-invariant forms of the 3D calculi.
In fact these commutation relations can now be obtained by 
differentiating the relations~(\ref{thfmcrs}) 
using~(\ref{hodif}). For the 4D calculi we 
obtain the following four sets of equations,
\begin{eqnarray}
\lefteqn{(\mathscr{C}_{st,j}A_{ij}-\mathscr{C}_{jk,i}(A_{js}A_{kt}+
B_{js}C_{kt}))\t_{s}\wedge\t_{t}=}\nonumber\;\;\;\;\;\;\\
& & 
A_{ij}(\t_{1}\wedge\t_{j}+\t_{j}\wedge\t_{1})+B_{ij}\t_{j}\wedge\t_{3}+C_{ij}\t_{2}\wedge\t_{j},
\nonumber\\
\lefteqn{(\mathscr{C}_{st,j}C_{ij}-\mathscr{C}_{jk,i}(C_{js}A_{kt}+
D_{js}C_{kt}))\t_{s}\wedge\t_{t}=}\nonumber\;\;\;\;\;\;\\
& & C_{ij}(\t_{4}\wedge\t_{j}+\t_{j}\wedge\t_{1})+D_{ij}\t_{j}\wedge\t_{3}+A_{ij}\t_{3}\wedge\t_{j},
\nonumber\\
\lefteqn{(\mathscr{C}_{st,j}B_{ij}-\mathscr{C}_{jk,i}(A_{js}B_{kt}+
B_{js}D_{kt}))\t_{s}\wedge\t_{t}=}\nonumber\;\;\;\;\;\;\\
& & B_{ij}(\t_{1}\wedge\t_{j}+\t_{j}\wedge\t_{4})+A_{ij}\t_{j}\wedge\t_{2}+D_{ij}\t_{2}\wedge\t_{j},
\nonumber\\
\lefteqn{(\mathscr{C}_{st,j}D_{ij}-\mathscr{C}_{jk,i}(C_{js}B_{kt}+
D_{js}D_{kt}))\t_{s}\wedge\t_{t}=}\nonumber\;\;\;\;\;\;\\
& & D_{ij}(\t_{4}\wedge\t_{j}+\t_{j}\wedge\t_{4})+C_{ij}\t_{j}\wedge\t_{2}+B_{ij}\t_{3}\wedge\t_{j},
\label{4Dththcrs}
\end{eqnarray}
where repeated indices are summed from 1 to 4.
In the 3D case the commutation relations are given by,
\begin{eqnarray}
\lefteqn{(\mathscr{C}_{st,j}A_{ij}-\mathscr{C}_{jk,i}(A_{js}A_{kt}+
B_{js}C_{kt}))\t_{s}\wedge\t_{t}=}\nonumber\;\;\;\;\;\;\\
& & A_{ij}(\t_{1}\wedge\t_{j}+\t_{j}\wedge\t_{1})+
B_{ij}\t_{j}\wedge\t_{3}+C_{ij}\t_{2}\wedge\t_{j},
\nonumber\\
\lefteqn{(\mathscr{C}_{st,j}C_{ij}-\mathscr{C}_{jk,i}(C_{js}A_{kt}+
D_{js}C_{kt}))\t_{s}\wedge\t_{t}=}\nonumber\;\;\;\;\;\;\\
& & C_{ij}(\t_{j}\wedge\t_{1}-\t_{1}\wedge\t_{j}-
2h\t_{3}\wedge\t_{j})+D_{ij}\t_{j}\wedge\t_{3}+A_{ij}\t_{3}\wedge\t_{j},
\nonumber\\
\lefteqn{(\mathscr{C}_{st,j}B_{ij}-\mathscr{C}_{jk,i}(A_{js}B_{kt}+
B_{js}D_{kt}))\t_{s}\wedge\t_{t}=}\nonumber\;\;\;\;\;\;\\
& & B_{ij}(\t_{1}\wedge\t_{j}-\t_{j}\wedge\t_{1}-
2h\t_{j}\wedge\t_{3})+A_{ij}\t_{j}\wedge\t_{2}+D_{ij}\t_{2}\wedge\t_{j},
\nonumber\\
\lefteqn{(\mathscr{C}_{st,j}D_{ij}-\mathscr{C}_{jk,i}(C_{js}B_{kt}+
D_{js}D_{kt})\t_{s}\wedge\t_{t}=}\nonumber\;\;\;\;\;\;\\
& & -D_{ij}(\t_{1}\wedge\t_{j}+\t_{j}\wedge\t_{1}+
2h\t_{3}\wedge\t_{j}+2h\t_{j}\wedge\t_{3})+C_{ij}\t_{j}\wedge\t_{2}+B_{ij}\t_{3}\wedge\t_{j},
\end{eqnarray}
where now repeated indices are summed from 1 to 3.

Recalling that $\Tr_{h}\Theta$ was the single biinvariant form in the 4D 
calculi, and also the defining characteristic of Woronowicz's 
bimodule map $\L$,~(\ref{defbraid}), we obtain a general relation 
which must be consistent with commutation relations between 1-forms 
in the 4D calculi,
\begin{eqnarray}
0=\Tr_{h}\Theta\wedge\Tr_{h}\Theta&=&\t_{1}\wedge\t_{1}+4h^{2}\t_{3}\wedge\t_{3}+\t_{4}\wedge\t_{4}+\nonumber\\
& &2h(\t_{1}\wedge\t_{3}+\t_{3}\wedge\t_{1})+\t_{1}\wedge\t_{4}+\t_{4}\wedge\t_{1}+
2h(\t_{3}\wedge\t_{4}+\t_{4}\wedge\t_{3}).\label{extra1}
\end{eqnarray}
Further, in any calculi where this bi-invariant form implements the 
differential in the sense of Remark~(\ref{reminner}), so that on 
arbitrary
1-forms $\w$ we have,
\begin{equation}
\dif \w=\frac{1}{\k}[Tr_{h}\Theta,\w]_{+},
\end{equation}
where $\k$ is some constant, we will have further relations,
\begin{equation}
[Tr_{h}\Theta,\t_{i}]=\k\sum_{j,k=1\ldots 
4}\mathscr{C}_{jk,i}\t_{j}\wedge\t_{k}.\label{extra2}
\end{equation}
Again these must be consistent with relations coming from~(\ref{4Dththcrs}).

The commutator expression for the quantum Lie bracket,~(\ref{rhs2}), requires 
that we know the explicit form of the matrix $\L$ whose components 
are given in~(\ref{braidcomp}) as $\L_{ij,st}=f_{it}(v_{sk})$. But we 
know the algebraic elements of the matrix $||v_{ij}||$ and therefore 
their expressions in terms of the $ABCD$ representation. This is all 
that is required. 

Finally, the structure constants, $\mathcal{C}_{ij,k}$, for the quantum Lie 
bracket, (\ref{rhs1}), now follow immediately from~(\ref{qstrcon}) 
as we know the $\mathscr{C}_{ij,k}$s and $\L_{ij,st}$s.
\label{MHProc}

\section{4-dimensional bicovariant calculi on $GL_{h,g}(2)$}
We summarise the result of applying the procedure of the previous 
section to $GL_{h,g}(2)$ in the following theorems.

\bth
There are three 1-parameter families of 4-dimensional first order 
bicovariant differential calculi on $GL_{h,g}(2)$ whose bimodules of forms 
are generated 
as left $GL_{h,g}(2)$-modules by the differentials of the quantum 
group generators. We will denote the three families by 
$\G^{\mathrm{4D}}_{1}$, 
$\G^{\mathrm{4D}}_{2}$ and $\G^{\mathrm{4D}}_{3}$. They are completely characterised by 
their respective $ABCD$ matrices.
\begin{description}
\item[$\G^{\mathrm{4D}}_{1}$]
The $ABCD$ matrices are given by,
\begin{align}
A&=\begin{pmatrix}
\frac{3z+2}{2} & 0 & \frac{-(3h+g)z-2h}{2} & \frac{-z}{2}\\
h(z+1) & z+1 & -h^{2}(z+1) & -h(z+1)\\
0 & 0 & z+1 & 0\\
\frac{z}{2} & 0 & \frac{(h-g)z+2h}{2} & \frac{z+2}{2}
\end{pmatrix},\nonumber\\
B&=\begin{pmatrix}
0 & z & h(h+g)(z+1) & 0\\
0 & (h+g)(z+1) & -hg(h+g)(z+1) & 0\\
0 & 0 & -(h+g)(z+1) & 0\\
0 & z & g(h+g)(z+1) & 0
\end{pmatrix},\nonumber\\
C&=\begin{pmatrix}
0 & 0 & z & 0\\
0 & 0 & 0 & 0\\
0 & 0 & 0 & 0\\
0 & 0 & z & 0
\end{pmatrix},\nonumber\\
D&=\begin{pmatrix}
\frac{z+2}{2} & 0 & \frac{-(h-g)z+2g}{2} & \frac{z}{2}\\
-g(z+1) & z+1 & -g^{2}(z+1) & g(z+1)\\
0 & 0 & z+1 & 0\\
\frac{-z}{2} & 0 & \frac{-(h+3g)z-2g}{2} & \frac{3z+2}{2}
\end{pmatrix}.
\end{align}
Here we must have $z\neq-1$ to ensure the invertibility of $\DD$. With 
$g=h$, the 
quantum determinant 
is central in the differential calculus for parameter values $z=0$ and $z=-2$. The 
differential of the quantum determinant is
\begin{equation}
\dif\mathcal{D}=\frac{z+2}{2}\mathcal{D}\Tr_{h}\Theta,
\end{equation}
and for $z\neq 0$ the calculi are inner,
\begin{align}
\dif a&=\frac{1}{2z}[\Tr_{h}\Theta,a],&\dif b&=\frac{1}{2z}[\Tr_{h}\Theta,b],\nonumber\\
\dif c&=\frac{1}{2z}[\Tr_{h}\Theta,c],&\dif d&=\frac{1}{2z}[\Tr_{h}\Theta,d].
\end{align}
\item[$\G^{\mathrm{4D}}_{2}$]
The $ABCD$ matrices are given by,
\begin{align}
A&=\begin{pmatrix}
\frac{z+2}{2} & 0 & \frac{(h+g)z-2h}{2} & \frac{z}{2}\\
h(z+1) & 1 & h((h+g)z-h) & h(z-1)\\
0 & 0 & 1 & 0\\
\frac{-z}{2} & 0 & \frac{-(h+g)z+2h}{2} & \frac{-z+2}{2}
\end{pmatrix},\nonumber\\
B&=\begin{pmatrix}
-hz & 0 & h(h+g)(1-z) & -hz\\
hgz & (h+g) & hg(h+g)(z-1) & ghz\\
z & 0 & (h+g)(z-1) & z\\
-gz & 0 & g(h+g)(1-z) & -gz
\end{pmatrix},\nonumber\\
C&=\begin{pmatrix}
0 & 0 & 0 & 0\\
z & 0 & (h+g)z & z\\
0 & 0 & 0 & 0\\
0 & 0 & 0 & 0
\end{pmatrix},\nonumber\\
D&=\begin{pmatrix}
\frac{-z+2}{2} & 0 & \frac{-(h+g)z+2g}{2} & \frac{-z}{2}\\
g(z-1) & 1 & g((h+g)z-g) & g(z+1)\\
0 & 0 & 1 & 0\\
\frac{z}{2} & 0 & \frac{(h+g)z-2g}{2} & \frac{z+2}{2}
\end{pmatrix}.
\end{align}
There is no restriction on the value of $z$ in this case. With $g=h$, 
the quantum 
determinant is central in the differential calculi for all values of 
$z$. The differential of the quantum determinant is
\begin{equation}
\dif\mathcal{D}=\frac{2-3z}{2}\mathcal{D}\Tr_{h}\Theta,
\end{equation}
but the calculi are not inner,
\begin{align}
[\Tr_{h}\Theta,a]&=0,&[\Tr_{h}\Theta,b]&=0,\nonumber\\
[\Tr_{h}\Theta,c]&=0,&[\Tr_{h}\Theta,d]&=0.
\end{align}
\item[$\G^{\mathrm{4D}}_{3}$]
The $ABCD$ matrices are given by,
\begin{align}
A&=\begin{pmatrix}
\frac{z+2}{2} & 0 & \frac{(h+g)z-2h}{2} & \frac{z}{2}\\
h & 1 & -h^{2} & -h\\
0 & 0 & 1 & 0\\
\frac{z}{2} & 0 & \frac{(h+g)z+2h}{2} & \frac{z+2}{2}
\end{pmatrix},&
B&=\begin{pmatrix}
0 & 0 & h(h+g) & 0\\
0 & (h+g) & -hg(h+g) & 0\\
0 & 0 & -(h+g) & 0\\
0 & 0 & g(h+g) & 0
\end{pmatrix},\nonumber\\
C&=\begin{pmatrix}
0 & 0 & 0 & 0\\
0 & 0 & 0 & 0\\
0 & 0 & 0 & 0\\
0 & 0 & 0 & 0
\end{pmatrix},&
D&=\begin{pmatrix}
\frac{z+2}{2} & 0 & \frac{(h+g)z+2g}{2} & \frac{z}{2}\\
-g & 1 & -g^{2} & g\\
0 & 0 & 1 & 0\\
\frac{z}{2} & 0 & \frac{(h+g)z-2g}{2} & \frac{z+2}{2}
\end{pmatrix}.
\end{align}
Like $\G^{\mathrm{4D}}_{1}$ we must have $z\neq-1$ to ensure the invertibility of $\DD$. 
With $g=h$, the 
quantum determinant 
is central in the differential calculus for parameter values $z=0$ and $z=-2$. The 
differential of the quantum determinant is
\begin{equation}
\dif\mathcal{D}=\frac{z+2}{2}\mathcal{D}\Tr_{h}\Theta,
\end{equation}
and again the calculi are not inner,
\begin{align}
[\Tr_{h}\Theta,a]&=za\Tr_{h}\Theta,&[\Tr_{h}\Theta,b]&=zb\Tr_{h}\Theta,\nonumber\\
[\Tr_{h}\Theta,c]&=zc\Tr_{h}\Theta,&[\Tr_{h}\Theta,d]&=zd\Tr_{h}\Theta.
\end{align}
\end{description}
For $z=0$, 
$\G^{\mathrm{4D}}_{1}=\G^{\mathrm{4D}}_{2}=\G^{\mathrm{4D}}_{3}$, 
whilst for all other parameter values the calculi are distinct.
\label{glabcd}
\ethe
\begin{proof}
The matrices were obtained by systematically extracting the implications 
of the Constraints~1--3. This laborious task is made possible by using 
the computer algebra package REDUCE~\cite{Redman}. REDUCE was also 
used to check the other results. The results regarding the centrality 
of the quantum determinant were obtained by investigating under what 
conditions $\DD=I$, since this is precisely the requirement imposed 
by the relations~(\ref{thfmcrs}). To obtain the expressions for 
the differential of the quantum determinant in the three calculi,
we differentiate the expression $\mathcal{D}=ad-bc+hac$, replace 
on the right hand side differentials by $\t_{i}$s 
through~(\ref{diftoth}), 
use the commutations relations between the $\t_{i}$ and the 
algebra generators provided by the $ABCD$ matrices to commute the 
left-invariant forms $\t_{i}$ to the right and finally 
straighten the algebra coefficients of the $\t_{i}$ to obtain the 
quoted results. 
\end{proof}

\bre
M\"{u}ller-Hoissen~\cite{MH1} obtained a corresponding result for first order 
differential calculi on the standard 2-parameter quantum group, 
$GL_{q,p}(2)$. In that case, there is just a \emph{single} 1-parameter 
family of calculi, and \emph{every} member of this family is inner.
\ere

Armed with the commutation relations between left-invariant forms and
generators provided by the classification 
of $ABCD$ matrices, we may now obtain commutation relations between 
the differentials of the generators and the generators. For example, 
in $\G^{\mathrm{4D}}_{1}$, starting with $\dif a a$, we replace $\dif 
a$ by it's expression in terms of left-invariant forms, commute these 
to the right and then replace the left-invariant forms by their 
expressions in terms of differentials to obtain,
\begin{align}
\dif a a=&(a+hc)\dif a+\frac{z}{4}\{9a+(12h-7g)c+\nonumber\\
&\mathcal{D}^{-1}\{bac-3a^{2}d+2(h-g)\{4bc^{2}-6adc-(h-2g)dc^{2}+3h(h-2g)c^{3}\}\}\}\dif a+
\nonumber\\
&\frac{z}{2}\mathcal{D}^{-1}\{a^{2}c+(h-g)\{2ac^{2}+(h-2g)c^{3}\}\}\dif 
b+\nonumber\\
&(ghc-ha)\dif c+\frac{z}{4}\{(7g-6h)a-((h-7g)(4h-3g)+gh)c+\nonumber\\
&\mathcal{D}^{-1}\{2ba^{2}-6ha^{3}+3(2h-3g)a^{2}d-(2h-3g)bac+\nonumber\\
&2(h-g)\{3(h-6g)adc-2(h-6g)bc^{2}-3g(h-2g)dc^{2}+9gh(h-2g)c^{3}\}\}\}\dif 
c\nonumber\\
&\frac{-z}{2}\mathcal{D}^{-1}\{a^{3}+(2h-3g)a^{2}c+(h-g)\{(h-6g)ac^{2}-3g(h-2g)c^{3}\}\}\dif d.
\end{align}
Most of the other commutation relations are much more complicated so 
we have chosen not to reproduce them here. A 
feature of this commutation relation, which is shared with the other 
differential-generator commutation relations for 
$\G^{\mathrm{4D}}_{1}$, $\G^{\mathrm{4D}}_{2}$ and 
$\G^{\mathrm{4D}}_{3}$, is that it is not 
quadratic, but that when $z=0$, in which case 
$\G^{\mathrm{4D}}_{1}=\G^{\mathrm{4D}}_{2}=\G^{\mathrm{4D}}_{3}$,
it simplifies drastically and does indeed become quadratic. Moreover, 
when $z=0$ these quadratic differential-generator commutation 
relations may be written in $R$-matrix form as\footnote{Of 
course, $\hat{R}^{-1}=\hat{R}$, 
but we write it this way to make 
comparison with the general results of other authors explicit.}
\begin{equation}
\hat{R}^{-1}\dif T_{1}T_{2}=T_{1}\dif T_{2}\hat{R}. \label{rmatdif}
\end{equation}
This $R$-matrix expression first appeared in the works of 
Schirrmacher~\cite{Schir} and Sudbery~\cite{Sud2} which were 
developments on the works of Manin~\cite{Man1,Man2} and 
Maltsiniotis~\cite{Malt}. These authors treated the co-plane, in our 
case $\AA^{2|0}_{-1}$, as the algebra of differentials $\dif x_{i}$ of the 
`coordinates' $x_{i}$ whose algebra is that of the plane, 
$\AA^{2|0}_{1}$. They then sought differential calculi expressed in 
terms of generators $T_{ij}$ and their differentials $\dif T_{ij}$ 
such that the plane and co-plane are invariant under the 
transformations $x_{i}\mapsto\sum T_{ij}x_{j}$ and $\dif 
x_{i}\mapsto\sum T_{ij}\dif x_{j}+\sum\dif T_{ij}x_{j}$. In the 
context of Jordanian quantum groups, in~\cite{Kar1} the author 
\emph{postulated}~(\ref{rmatdif}) as 
the relation defining commutation relations between differentials and 
generators for $SL_{h}(2)$.

\bre
In~\cite{MHR}, where the differential calculus on the standard quantum 
group $GL_{q,p}(2)$ is considered, the authors also observe that 
the differential-generator commutation relations are not quadratic 
for general calculi in the 1-parameter family, but that once their free 
parameter is fixed to zero, quadratic relations are obtained. Moreover 
in this case they also recover commutation relations with an 
$R$-matrix expression.
\ere

\bre
From~(\ref{rmatdif}) it is a simple matter to demonstrate that for 
the $z=0$ calculus on $GL_{h,g}(2)$, the commutation relations between 
the left-invariant forms and the quantum group generators 
may be written as,
\begin{equation}
\Theta_{1}T_{2}=T_{2}R_{21}\Theta_{1}R_{12},
\end{equation}
where
\begin{equation}
\Theta=\begin{pmatrix}
\t_{1} & \t_{2}\\
\t_{3} & \t_{4}
\end{pmatrix}.
\end{equation}
Incidentally, this is precisely the relation we obtain if we attempt a 
naive application of Jurco's construction with, in the notation 
of~\cite{Jur2}, $\G=\G_{1}^{c}\ot_{A}\G_{1}$.
\ere

Turning now to the higher order calculi, we have the following result 
describing commutation relations between the left-invariant forms.

\bth
The commutation relations between the left-invariant forms in the 
external bicovariant graded algebras $\Omega^{\mathrm{4D}}_{1}$, 
$\Omega^{\mathrm{4D}}_{2}$ and $\Omega^{\mathrm{4D}}_{3}$ built 
respectively upon the three families of 
first order calculi, $\G^{\mathrm{4D}}_{1}$, 
$\G^{\mathrm{4D}}_{2}$ and $\G^{\mathrm{4D}}_{3}$ using Woronowicz's 
theory as described in Section~\ref{WorThr} are as follows,
\begin{description}
\item[$\Omega^{\mathrm{4D}}_{1}$]
\begin{align}
\t_{3}\wedge\t_{3}&=0\nonumber\\
\t_{3}\wedge\t_{4}&=-\frac{2z+1}{z+1}\t_{4}\wedge\t_{3}+\frac{z}{z+1}\t_{1}\wedge\t_{3},
\nonumber\\
\t_{3}\wedge\t_{1}&=-\frac{1}{z+1}\t_{1}\wedge\t_{3}-\frac{z}{z+1}\t_{4}\wedge\t_{3},
\nonumber\\
\t_{3}\wedge\t_{2}&=-\t_{2}\wedge\t_{3}+(h+g)\t_{1}\wedge\t_{3}-(h+g)\t_{4}\wedge\t_{3},
\nonumber\\
\t_{4}\wedge\t_{4}&=\frac{z}{z+1}\t_{2}\wedge\t_{3}-\frac{z(h+g)}{z+1}\t_{1}\wedge\t_{3}+
\frac{z(h+g)}{z+1}\t_{4}\wedge\t_{3},\nonumber\\
\t_{4}\wedge\t_{1}&=-\t_{1}\wedge\t_{4}-\frac{z(h+g)}{z+1}\t_{1}\wedge\t_{3}+
\frac{z(h+g)}{z+1}\t_{4}\wedge\t_{3},\nonumber\\
\t_{4}\wedge\t_{2}&=-\frac{2z+1}{z+1}\t_{2}\wedge\t_{4}+\frac{z}{z+1}\t_{2}\wedge\t_{1}+
\frac{(2z+1)(h+g)}{z+1}\t_{2}\wedge\t_{3}\nonumber\\
&\quad-(h+g)^{2}\t_{1}\wedge\t_{3}+
(h+g)^{2}\t_{4}\wedge\t_{3},\nonumber\\
\t_{1}\wedge\t_{1}&=-\frac{z}{z+1}\t_{2}\wedge\t_{3},\nonumber\\
\t_{1}\wedge\t_{2}&=-\frac{1}{z+1}\t_{2}\wedge\t_{1}-\frac{z}{z+1}\t_{2}\wedge\t_{4}-
\frac{(h+g)}{z+1}\t_{2}\wedge\t_{3},\nonumber\\
\t_{2}\wedge\t_{2}&=(h+g)\t_{1}\wedge\t_{2}-(h+g)\t_{2}\wedge\t_{4}+
(h+g)^{2}\t_{2}\wedge\t_{3},
\end{align}
\item[$\Omega^{\mathrm{4D}}_{2}$]
\begin{align}
\t_{3}\wedge\t_{3}&=0\nonumber\\
\t_{3}\wedge\t_{4}&=-\t_{4}\wedge\t_{3},
\nonumber\\
\t_{3}\wedge\t_{1}&=-\t_{1}\wedge\t_{3},
\nonumber\\
\t_{3}\wedge\t_{2}&=-\t_{2}\wedge\t_{3}+(h+g)\t_{1}\wedge\t_{3}-(h+g)\t_{4}\wedge\t_{3},
\nonumber\\
\t_{4}\wedge\t_{4}&=0,\nonumber\\
\t_{4}\wedge\t_{1}&=-\t_{1}\wedge\t_{4},\nonumber\\
\t_{4}\wedge\t_{2}&=-\t_{2}\wedge\t_{4}+
(h+g)\t_{2}\wedge\t_{3}-(h+g)^{2}\t_{1}\wedge\t_{3}+
(h+g)^{2}\t_{4}\wedge\t_{3},\nonumber\\
\t_{1}\wedge\t_{1}&=0,\nonumber\\
\t_{1}\wedge\t_{2}&=-\t_{2}\wedge\t_{1}-
(h+g)\t_{2}\wedge\t_{3},\nonumber\\
\t_{2}\wedge\t_{2}&=(h+g)\t_{1}\wedge\t_{2}-(h+g)\t_{2}\wedge\t_{4}+
(h+g)^{2}\t_{2}\wedge\t_{3},
\end{align}
\item[$\Omega^{\mathrm{4D}}_{3}$]
\begin{align}
\t_{3}\wedge\t_{3}&=0\nonumber\\
\t_{3}\wedge\t_{4}&=-\t_{4}\wedge\t_{3},
\nonumber\\
\t_{3}\wedge\t_{1}&=-\t_{1}\wedge\t_{3},
\nonumber\\
\t_{3}\wedge\t_{2}&=-\t_{2}\wedge\t_{3}+(h+g)\t_{1}\wedge\t_{3}-(h+g)\t_{4}\wedge\t_{3},
\nonumber\\
\t_{4}\wedge\t_{4}&=0,\nonumber\\
\t_{4}\wedge\t_{1}&=-\t_{1}\wedge\t_{4},\nonumber\\
\t_{4}\wedge\t_{2}&=-\t_{2}\wedge\t_{4}+
(h+g)\t_{2}\wedge\t_{3}-(h+g)^{2}\t_{1}\wedge\t_{3}+
(h+g)^{2}\t_{4}\wedge\t_{3},\nonumber\\
\t_{1}\wedge\t_{1}&=0,\nonumber\\
\t_{1}\wedge\t_{2}&=-\t_{2}\wedge\t_{1}-
(h+g)\t_{2}\wedge\t_{3},\nonumber\\
\t_{2}\wedge\t_{2}&=(h+g)\t_{1}\wedge\t_{2}-(h+g)\t_{2}\wedge\t_{4}+
(h+g)^{2}\t_{2}\wedge\t_{3}.
\end{align}
\end{description}
The relations in each case are consistent with the 
respective relations~(\ref{thfmcrs}). They are also consistent with 
the relation~(\ref{extra1}) and in the case of 
$\Omega^{\mathrm{4D}}_{1}$, the relations~(\ref{extra2}). Further, they 
are such that 
$\{\t_{2}^{\a}\t_{1}^{\b}\t_{4}^{\g}\t_{3}^{\d}:\a,\b,\g,\d\in\{0,1\}\}$ 
is a basis for the exterior algebra of forms in each case.
\label{wedcoms}
\ethe
\begin{proof}
Once again we used REDUCE to check these results based on the 
discussion of Section~\ref{MHProc}. The 16 equations~(\ref{4Dththcrs}) 
were treated, in each of the cases, $\G^{\mathrm{4D}}_{1}$, 
$\G^{\mathrm{4D}}_{2}$ and $\G^{\mathrm{4D}}_{3}$ linear relations 
between the 10 `mis-ordered' 2-forms, $\t_{3}\wedge\t_{3}$, 
$\t_{3}\wedge\t_{4}$, $\t_{3}\wedge\t_{1}$, $\t_{3}\wedge\t_{2}$, 
$\t_{4}\wedge\t_{4}$, $\t_{4}\wedge\t_{1}$, $\t_{4}\wedge\t_{2}$, 
$\t_{1}\wedge\t_{1}$, $\t_{1}\wedge\t_{2}$ and $\t_{2}\wedge\t_{2}$, 
and the 6 `ordered' 2-forms, $\t_{2}\wedge\t_{1}$, 
$\t_{2}\wedge\t_{4}$, $\t_{2}\wedge\t_{3}$, $\t_{1}\wedge\t_{4}$, 
$\t_{1}\wedge\t_{3}$ and $\t_{4}\wedge\t_{3}$. These linear relations 
were then solved for the 10 mis-ordered 2-forms yielding in each case 
the single solution presented. Consistency with the 
relations~(\ref{thfmcrs}) in each of the three cases was checked by commuting 
the generators through the relations and observing that no further 
conditions were incurred. The other consistency conditions were again 
checked by direct computation. Finally, observing that the relations 
are compatible with the ordering 
$\t_{2}\prec\t_{1}\prec\t_{4}\prec\t_{3}$, we may use the 
Diamond Lemma to prove the statement about the bases.
\end{proof}

\bre
It is interesting to note here that the relations in 
$\G^{\mathrm{4D}}_{2}$ and $\G^{\mathrm{4D}}_{3}$ are the same and 
indeed could be obtained from those in $\G^{\mathrm{4D}}_{1}$ by 
setting $z=0$. 
\ere

\bre
In contrast with the results in Theorem~\ref{wedcoms}, in the work of 
M\"{u}ller-Hoissen and Reuten on $GL_{q,p}(2)$~\cite{MHR} 
the corresponding commutation relations between left-invariant forms 
exhibited ordering circles which introduced further constraints on 
the free parameter of their family of calculi.
\ere

Let us now describe the quantum Lie brackets in the quantum Lie 
algebras $\mathcal{L}^{\mathrm{4D}}_{1}$, 
$\mathcal{L}^{\mathrm{4D}}_{2}$ and $\mathcal{L}^{\mathrm{4D}}_{3}$ 
dual to the bicovariant bimodules $\G^{\mathrm{4D}}_{1}$, 
$\G^{\mathrm{4D}}_{2}$ and $\G^{\mathrm{4D}}_{3}$ respectively. 

\bth
The quantum Lie brackets and quantum commutators for the 
quantum Lie algebras $\mathcal{L}^{\mathrm{4D}}_{1}$, 
$\mathcal{L}^{\mathrm{4D}}_{2}$ and $\mathcal{L}^{\mathrm{4D}}_{3}$ 
as described in Section~\ref{WorThr} are as follows:
\begin{description}
\item[$\mathcal{L}^{\mathrm{4D}}_{1}$]
The bracket relations are,
\begin{align}
[\chi_{1},\chi_{1}]&=0, \nonumber \\
[\chi_{1},\chi_{2}]&=\frac{1}{z+1}\chi_{2}, \nonumber \\
[\chi_{1},\chi_{3}]&=-\frac{1}{z+1}\chi_{3}+\frac{h+g}{z+1}\chi_{4},\nonumber\\
[\chi_{1},\chi_{4}]&=0,\nonumber\\
[\chi_{2},\chi_{1}]&=-\frac{1}{z+1}\chi_{2},\nonumber\\
[\chi_{2},\chi_{2}]&=0,\nonumber\\
[\chi_{2},\chi_{3}]&=\frac{1}{z+1}\chi_{1}-\frac{h+g}{z+1}\chi_{2}-\frac{1}{z+1}\chi_{4},\nonumber\\
[\chi_{2},\chi_{4}]&=\frac{1}{z+1}\chi_{2},\nonumber\\
[\chi_{3},\chi_{1}]&=-\frac{h+g}{z+1}\chi_{1}+\frac{1}{z+1}\chi_{3},\nonumber\\
[\chi_{3},\chi_{2}]&=-\frac{1}{z+1}\chi_{1}-\frac{h+g}{z+1}\chi_{2}+\frac{1}{z+1}\chi_{4},\nonumber\\
[\chi_{3},\chi_{3}]&=\frac{(h+g)^{2}}{z+1}\chi_{1}-2\frac{h+g}{z+1}\chi_{3}+
\frac{(h+g)^{2}}{z+1}\chi_{4},\nonumber\\
[\chi_{3},\chi_{4}]&=\frac{h+g}{z+1}\chi_{1}-\frac{1}{z+1}\chi_{3},\nonumber\\
[\chi_{4},\chi_{1}]&=0,\nonumber\\
[\chi_{4},\chi_{2}]&=-\frac{1}{z+1}\chi_{2},\nonumber\\
[\chi_{4},\chi_{3}]&=\frac{1}{z+1}\chi_{3}-\frac{h+g}{z+1}\chi_{4},\nonumber\\
[\chi_{4},\chi_{4}]&=0,
\end{align}
and the commutators are
\begin{align*}
[\chi_{1},\chi_{1}]&=0,\nonumber\\
[\chi_{1},\chi_{2}]&=\chi_{1}\chi_{2}-\Bigl(\frac{z}{z+1}\chi_{1}\chi_{2}+\chi_{2}\chi_{1}
+(h+g)\chi_{2}\chi_{2}+\frac{z}{z+1}\chi_{4}\chi_{2}\Bigr),\nonumber\\
[\chi_{1},\chi_{3}]&=\chi_{1}\chi_{3}-\Bigl(\frac{z}{z+1}\chi_{1}\chi_{3}+
\frac{z(h+g)}{z+1}\chi_{1}\chi_{4}+\chi_{3}\chi_{1}-(h+g)\chi_{3}\chi_{2}+(h+g)^{2}\chi_{4}\chi_{2}
\nonumber\\
&\quad-\frac{z}{z+1}\chi_{4}\chi_{3}+\frac{z(h+g)}{z+1}\chi_{4}\chi_{4}\Bigr),\nonumber\\
[\chi_{1},\chi_{4}]&=\chi_{1}\chi_{4}-\chi_{4}\chi_{1},\nonumber\\
[\chi_{2},\chi_{1}]&=\chi_{2}\chi_{1}-\Bigl(\frac{1}{z+1}\chi_{1}\chi_{2}-(h+g)\chi_{2}\chi_{2}
-\frac{z}{z+1}\chi_{4}\chi_{2}\Bigr),\nonumber\\
[\chi_{2},\chi_{2}]&=0,\nonumber\\
\end{align*}
\begin{align}
[\chi_{2},\chi_{3}]&=\chi_{2}\chi_{3}-\Bigl(\frac{z}{z+1}\chi_{1}\chi_{1}+
\frac{h+g}{z+1}\chi_{1}\chi_{2}-\frac{z}{z+1}\chi_{1}\chi_{4}-(h+g)^{2}\chi_{2}\chi_{2}+
\chi_{3}\chi_{2}\nonumber\\
&\quad+\frac{z}{z+1}\chi_{4}\chi_{1}-\frac{(h+g)(2z+1)}{z+1}\chi_{4}\chi_{2}-
\frac{z}{z+1}\chi_{4}\chi_{4}\Bigr),\nonumber\\
[\chi_{2},\chi_{4}]&=\chi_{2}\chi_{4}-\Bigl(\frac{z}{z+1}\chi_{1}\chi_{2}+(h+g)\chi_{2}\chi_{2}+
\frac{2z+1}{z+1}\chi_{4}\chi_{2}\Bigr),\nonumber\\
[\chi_{3},\chi_{1}]&=\chi_{3}\chi_{1}-\Bigl(-\frac{z(h+g)}{z+1}\chi_{1}\chi_{1}+
\frac{2z+1}{z+1}\chi_{1}\chi_{3}-(h+g)^{2}\chi_{2}\chi_{1}+(h+g)\chi_{2}\chi_{3}\nonumber\\
&\quad-\frac{z(h+g)}{z+1}\chi_{4}\chi_{1}+\frac{z}{z+1}\chi_{4}\chi_{3}\Bigr),\nonumber\\
[\chi_{3},\chi_{2}]&=\chi_{3}\chi_{2}-\Bigl(-\frac{z}{z+1}\chi_{1}\chi_{1}-
\frac{z(h+g)}{z+1}\chi_{1}\chi_{2}+\frac{z}{z+1}\chi_{1}\chi_{4}-(h+g)\chi_{2}\chi_{1}\nonumber\\
&\quad-(h+g)^{2}\chi_{2}\chi_{2}+\chi_{2}\chi_{3}+(h+g)\chi_{2}\chi_{4}-
\frac{z}{z+1}\chi_{4}\chi_{1}-
\frac{z(h+g)}{z+1}\chi_{4}\chi_{2}\nonumber\\
&\quad+\frac{z}{z+1}\chi_{4}\chi_{4}\Bigr),\nonumber\\
[\chi_{3},\chi_{3}]&=\chi_{3}\chi_{3}-\Bigl(\frac{z(h+g)^{2}}{z+1}\chi_{1}\chi_{1}-
\frac{(h+g)(3z+1)}{z+1}\chi_{1}\chi_{3}+\frac{(h+g)^{2}(2z+1)}{z+1}\chi_{1}\chi_{4}\nonumber\\
&\quad+(h+g)^{3}\chi_{2}\chi_{1}-(h+g)^{2}\chi_{2}\chi_{3}+(h+g)\chi_{3}\chi_{1}-
(h+g)^{2}\chi_{3}\chi_{2}+\chi_{3}\chi_{3}\nonumber\\
&\quad-(h+g)\chi_{3}\chi_{4}-\frac{(h+g)^{2}}{z+1}\chi_{4}\chi_{1}+(h+g)^{3}\chi_{4}\chi_{2}
-\frac{(h+g)(z-1)}{z+1}\chi_{4}\chi_{3}\nonumber\\
&\quad+\frac{z(h+g)^{2}}{z+1}\chi_{4}\chi_{4}\Bigr),\nonumber\\
[\chi_{3},\chi_{4}]&=\chi_{3}\chi_{4}-\Bigl(\frac{z(h+g)}{z+1}\chi_{1}\chi_{1}-
\frac{z}{z+1}\chi_{1}\chi_{3}+(h+g)^{2}\chi_{2}\chi_{1}-(h+g)\chi_{2}\chi_{3}\nonumber\\
&\quad+\frac{z(h+g)}{z+1}\chi_{4}\chi_{1}+\frac{1}{z+1}\chi_{4}\chi_{3}\Bigr),\nonumber\\
[\chi_{4},\chi_{1}]&=\chi_{4}\chi_{1}-\chi_{1}\chi_{4},\nonumber\\
[\chi_{4},\chi_{2}]&=\chi_{4}\chi_{2}-\Bigl(-\frac{z}{z+1}\chi_{1}\chi_{2}-
(h+g)\chi_{2}\chi_{2}+\chi_{2}\chi_{4}-\frac{z}{z+1}\chi_{4}\chi_{2}\Bigr),\nonumber\\
[\chi_{4},\chi_{3}]&=\chi_{4}\chi_{3}-\Bigl(\frac{z}{z+1}\chi_{1}\chi_{3}-
\frac{z(h+g)}{z+1}\chi_{1}\chi_{4}+(h+g)\chi_{3}\chi_{2}+\chi_{3}\chi_{4}-
(h+g)^{2}\chi_{4}\chi_{2}\nonumber\\
&\quad+\frac{z}{z+1}\chi_{4}\chi_{3}-\frac{z(h+g)}{z+1}\chi_{4}\chi_{4}\Bigr),\nonumber\\
[\chi_{4},\chi_{4}]&=0.
\end{align}
\item[$\mathcal{L}^{\mathrm{4D}}_{2}$]
The bracket relations are,
\begin{align*}
[\chi_{1},\chi_{1}]&=z\chi_{1}-z\chi_{4},\nonumber\\
[\chi_{1},\chi_{2}]&=\chi_{2},\nonumber\\
[\chi_{1},\chi_{3}]&=z(h+g)\chi_{1}-\chi_{3}-(z-1)(h+g)\chi_{4},\nonumber\\
[\chi_{1},\chi_{4}]&=z\chi_{1}-z\chi_{4},\nonumber\\
[\chi_{2},\chi_{1}]&=(2z-1)\chi_{2},\nonumber\\
[\chi_{2},\chi_{2}]&=0,\nonumber\\
[\chi_{2},\chi_{3}]&=\chi_{1}+(2z-1)(h+g)\chi_{2}-\chi_{4},\nonumber\\
[\chi_{2},\chi_{4}]&=(2z+1)\chi_{2},\nonumber\\
\end{align*}
\begin{align}
[\chi_{3},\chi_{1}]&=-(z+1)(h+g)\chi_{1}+(2z+1)\chi_{3}-z(h+g)\chi_{4},\nonumber\\
[\chi_{3},\chi_{2}]&=-\chi_{1}-(h+g)\chi_{2}+\chi_{4},\nonumber\\
[\chi_{3},\chi_{3}]&=-(z-1)(h+g)^{2}\chi_{1}+2(z-1)(h+g)\chi_{3}
-(z-1)(h+g)^{2}\chi_{4},\nonumber\\
[\chi_{3},\chi_{4}]&=-(z-1)(h+g)\chi_{1}+(2z-1)\chi_{3}-z(h+g)\chi_{4},\nonumber\\
[\chi_{4},\chi_{1}]&=-z\chi_{1}+z\chi_{4},\nonumber\\
[\chi_{4},\chi_{2}]&=-\chi_{2},\nonumber\\
[\chi_{4},\chi_{3}]&=-z(h+g)\chi_{1}+\chi_{3}+(z-1)(h+g)\chi_{4},\nonumber\\
[\chi_{4},\chi_{4}]&=-z\chi_{1}+z\chi_{4},
\end{align}
and the commutators are,
\begin{align*}
[\chi_{1},\chi_{1}]&=\chi_{1}\chi_{1}-\Bigl(\chi_{1}\chi_{1}+z(h+g)\chi_{2}\chi_{1}-
z\chi_{2}\chi_{3}+z\chi_{3}\chi_{2}-z(h+g)\chi_{4}\chi_{2}\Bigr),\nonumber\\
[\chi_{1},\chi_{2}]&=\chi_{1}\chi_{2}-\Bigl(\chi_{2}\chi_{1}+(h+g)\chi_{2}\chi_{2}\Bigr),\nonumber\\
[\chi_{1},\chi_{3}]&=\chi_{1}\chi_{3}-\Bigl(z(h+g)^{2}\chi_{2}\chi_{1}-z(h+g)\chi_{2}\chi_{3}+
\chi_{3}\chi_{1}+(h+g)(z-1)\chi_{3}\chi_{2}\nonumber\\
&\quad-(h+g)^{2}(z-1)\chi_{4}\chi_{2}\Bigr),\nonumber\\
[\chi_{1},\chi_{4}]&=\chi_{1}\chi_{4}-\Bigl(z(h+g)\chi_{2}\chi_{1}-z\chi_{2}\chi_{3}+
z\chi_{3}\chi_{2}+\chi_{4}\chi_{1}-z(h+g)\chi_{4}\chi_{2}\Bigr),\nonumber\\
[\chi_{2},\chi_{1}]&=\chi_{2}\chi_{1}-\Bigl(-(z-1)\chi_{1}\chi_{2}+z\chi_{2}\chi_{1}+
(h+g)(2z-1)\chi_{2}\chi_{2}-z\chi_{2}\chi_{4}\nonumber\\
&\quad+z\chi_{4}\chi_{2}\Bigr),\nonumber\\
[\chi_{2},\chi_{2}]&=0,\nonumber\\
[\chi_{2},\chi_{3}]&=\chi_{2}\chi_{3}-\Bigl(-(h+g)(z-1)\chi_{1}\chi_{2}+z(h+g)\chi_{2}\chi_{1}+
(h+g)^{2}(2z-1)\chi_{2}\chi_{2}\nonumber\\
&\quad-z(h+g)\chi_{2}\chi_{4}+\chi_{3}\chi_{2}+(h+g)(z-1)\chi_{4}\chi_{2}\Bigr),\nonumber\\
[\chi_{2},\chi_{4}]&=\chi_{2}\chi_{4}-\Bigl(-z\chi_{1}\chi_{2}+z\chi_{2}\chi_{1}+
(h+g)(2z+1)\chi_{2}\chi_{2}-z\chi_{2}\chi_{4}+(z+1)\chi_{4}\chi_{2}\Bigr),\nonumber\\
[\chi_{3},\chi_{1}]&=\chi_{3}\chi_{1}-\Bigl((z+1)\chi_{1}\chi_{3}-z(h+g)\chi_{1}\chi_{4}-
(h+g)^{2}(z+1)\chi_{2}\chi_{1}+(h+g)(z+1)\chi_{2}\chi_{3}\nonumber\\
&\quad-z\chi_{3}\chi_{1}+
z(h+g)\chi_{3}\chi_{2}+z\chi_{3}\chi_{4}+z(h+g)\chi_{4}\chi_{1}-
z(h+g)^{2}\chi_{4}\chi_{2}-z\chi_{4}\chi_{3}\Bigr),\nonumber\\
[\chi_{3},\chi_{2}]&=\chi_{3}\chi_{2}-\Bigl(-(h+g)\chi_{2}\chi_{1}-
(h+g)^{2}\chi_{2}\chi_{2}+
\chi_{2}\chi_{3}+(h+g)\chi_{2}\chi_{4}\Bigr),\nonumber\\
[\chi_{3},\chi_{3}]&=\chi_{3}\chi_{3}-\Bigl((h+g)(z-1)\chi_{1}\chi_{3}-
(h+g)^{2}(z-1)\chi_{1}\chi_{4}
-(h+g)^{3}(z-1)\chi_{2}\chi_{1}\nonumber\\
&\quad+(h+g)^{2}(z-1)\chi_{2}\chi_{3}-(h+g)(z-1)\chi_{3}\chi_{1}+
(h+g)^{2}(z-1)\chi_{3}\chi_{2}+\chi_{3}\chi_{3}\nonumber\\
&\quad+(h+g)(z-1)\chi_{3}\chi_{4}+
(h+g)^{2}(z-1)\chi_{4}\chi_{1}-(h+g)^{3}(z-1)\chi_{4}\chi_{2}\nonumber\\
&\quad-(h+g)(z-1)\chi_{4}\chi_{3}\Bigr),\nonumber\\
[\chi_{3},\chi_{4}]&=\chi_{3}\chi_{4}-\Bigl(z\chi_{1}\chi_{3}-z(h+g)\chi_{1}\chi_{4}-
(h+g)^{2}(z-1)\chi_{2}\chi_{1}+(h+g)(z-1)\chi_{2}\chi_{3}\nonumber\\
&\quad-z\chi_{3}\chi_{1}+z(h+g)\chi_{3}\chi_{2}+
z\chi_{3}\chi_{4}+z(h+g)\chi_{4}\chi_{1}-z(h+g)^{2}\chi_{4}\chi_{2}\nonumber\\
&\quad-(z-1)\chi_{4}\chi_{3}\Bigr),\nonumber\\
\end{align*}
\begin{align}
[\chi_{4},\chi_{1}]&=\chi_{4}\chi_{1}-\Bigl(\chi_{1}\chi_{4}-z(h+g)\chi_{2}\chi_{1}+z\chi_{2}\chi_{3}
-z\chi_{3}\chi_{2}+z(h+g)\chi_{4}\chi_{2}\Bigr),\nonumber\\
[\chi_{4},\chi_{2}]&=\chi_{4}\chi_{2}-\Bigl(-(h+g)\chi_{2}\chi_{2}+\chi_{2}\chi_{4}\Bigr),\nonumber\\
[\chi_{4},\chi_{3}]&=\chi_{4}\chi_{3}-\Bigl(-z(h+g)^{2}+z(h+g)\chi_{2}\chi_{3}-
(h+g)(z-1)\chi_{3}\chi_{2}+\chi_{3}\chi_{4}\nonumber\\
&\quad+(h+g)^{2}(z-1)\chi_{4}\chi_{2}\Bigr),\nonumber\\
[\chi_{4},\chi_{4}]&=\chi_{4}\chi_{4}-\Bigl(-z(h+g)\chi_{2}\chi_{1}+z\chi_{2}\chi_{3}-
z\chi_{3}\chi_{2}+z(h+g)\chi_{4}\chi_{2}+\chi_{4}\chi_{4}\Bigr).
\end{align}
\item[$\mathcal{L}^{\mathrm{4D}}_{3}$]
The bracket relations are,
\begin{align}
[\chi_{1},\chi_{1}]&=0,\nonumber\\
[\chi_{1},\chi_{2}]&=\chi_{2},\nonumber\\
[\chi_{1},\chi_{3}]&=-\chi_{3}+(h+g)\chi_{4},\nonumber\\
[\chi_{1},\chi_{4}]&=0,\nonumber\\
[\chi_{2},\chi_{1}]&=-\chi_{2},\nonumber\\
[\chi_{2},\chi_{2}]&=0,\nonumber\\
[\chi_{2},\chi_{3}]&=\chi_{1}-(h+g)\chi_{2}-\chi_{4},\nonumber\\
[\chi_{2},\chi_{4}]&=\chi_{2},\nonumber\\
[\chi_{3},\chi_{1}]&=-(h+g)\chi_{1}+\chi_{3},\nonumber\\
[\chi_{3},\chi_{2}]&=-\chi_{1}-(h+g)\chi_{2}+\chi_{4},\nonumber\\
[\chi_{3},\chi_{3}]&=(h+g)^{2}\chi_{1}-2(h+g)\chi_{3}+(h+g)^{2}\chi_{4},\nonumber\\
[\chi_{3},\chi_{4}]&=(h+g)\chi_{1}-\chi_{3},\nonumber\\
[\chi_{4},\chi_{1}]&=0,\nonumber\\
[\chi_{4},\chi_{2}]&=-\chi_{2},\nonumber\\
[\chi_{4},\chi_{3}]&=\chi_{3}-(h+g)\chi_{4},\nonumber\\
[\chi_{4},\chi_{4}]&=0,
\end{align}
and the commutators are,
\begin{align*}
[\chi_{1},\chi_{1}]&=0,\nonumber\\
[\chi_{1},\chi_{2}]&=\chi_{1}\chi_{2}-\Bigl(\chi_{2}\chi_{1}+(h+g)\chi_{2}\chi_{2}\Bigr),\nonumber\\
[\chi_{1},\chi_{3}]&=\chi_{1}\chi_{3}-\Bigl(\chi_{3}\chi_{1}-(h+g)\chi_{3}\chi_{2}+
(h+g)^{2}\chi_{4}\chi_{2}\Bigr),\nonumber\\
[\chi_{1},\chi_{4}]&=0,\nonumber\\
[\chi_{2},\chi_{1}]&=\chi_{2}\chi_{1}-\Bigl(\chi_{1}\chi_{2}-(h+g)\chi_{2}\chi_{2}\Bigr),\nonumber\\
[\chi_{2},\chi_{2}]&=0,\nonumber\\
[\chi_{2},\chi_{3}]&=\chi_{2}\chi_{3}-\Bigl((h+g)\chi_{1}\chi_{2}-(h+g)^{2}\chi_{2}\chi_{2}+
\chi_{3}\chi_{2}-(h+g)\chi_{4}\chi_{2}\Bigr),\nonumber\\
[\chi_{2},\chi_{4}]&=\chi_{2}\chi_{4}-\Bigl((h+g)\chi_{2}\chi_{2}+\chi_{4}\chi_{2}\Bigr),\nonumber\\
[\chi_{3},\chi_{1}]&=\chi_{3}\chi_{1}-\Bigl(\chi_{1}\chi_{3}-(h+g)^{2}\chi_{2}\chi_{1}+
(h+g)\chi_{2}\chi_{3}\Bigr),\nonumber\\
\end{align*}
\begin{align}
[\chi_{3},\chi_{2}]&=\chi_{3}\chi_{2}-\Bigl(-(h+g)\chi_{2}\chi_{1}-(h+g)^{2}\chi_{2}\chi_{2}+
\chi_{2}\chi_{3}+(h+g)\chi_{2}\chi_{4}\Bigr),\nonumber\\
[\chi_{3},\chi_{3}]&=\chi_{3}\chi_{3}-\Bigl(-(h+g)\chi_{1}\chi_{3}+(h+g)^{2}\chi_{1}\chi_{4}+
(h+g)^{3}\chi_{2}\chi_{1}-(h+g)^{2}\chi_{2}\chi_{3}\nonumber\\
&\quad+(h+g)\chi_{3}\chi_{1}-(h+g)^{2}\chi_{3}\chi_{2}+
\chi_{3}\chi_{3}-(h+g)\chi_{3}\chi_{4}-(h+g)^{2}\chi_{4}\chi_{1}\nonumber\\
&\quad+(h+g)^{3}\chi_{4}\chi_{2}+
(h+g)\chi_{4}\chi_{3}\Bigr),\nonumber\\
[\chi_{3},\chi_{4}]&=\chi_{3}\chi_{4}-\Bigl((h+g)^{2}\chi_{2}\chi_{1}-(h+g)\chi_{2}\chi_{3}+
\chi_{4}\chi_{3}\Bigr),\nonumber\\
[\chi_{4},\chi_{1}]&=0,\nonumber\\
[\chi_{4},\chi_{2}]&=\chi_{4}\chi_{2}-\Bigl(-(h+g)\chi_{2}\chi_{2}+\chi_{2}\chi_{4}\Bigr),\nonumber\\
[\chi_{4},\chi_{3}]&=\chi_{4}\chi_{3}-\Bigl((h+g)\chi_{3}\chi_{2}+chi_{3}\chi_{4}-
(h+g)^{2}\chi_{4}\chi_{2}\Bigr),\nonumber\\
[\chi_{4},\chi_{4}]&=0.
\end{align}
\end{description}
\ethe

The following result reveals that the relations in the 
universal enveloping algebras $U(\mathcal{L}^{\mathrm{4D}}_{1})$, 
$U(\mathcal{L}^{\mathrm{4D}}_{2})$ and $U(\mathcal{L}^{\mathrm{4D}}_{3})$ 
reflect the structure of the commutation relations of the 
left-invariant forms presented in Theorem~\ref{wedcoms}. Indeed the 
relations in $U(\mathcal{L}^{\mathrm{4D}}_{2})$ and 
$U(\mathcal{L}^{\mathrm{4D}}_{3})$ are identical and can be obtained 
from those in $U(\mathcal{L}^{\mathrm{4D}}_{1})$ by setting $z=0$

\bth
The relations in the universal enveloping algebras, $U(\mathcal{L}^{\mathrm{4D}}_{1})$, 
$U(\mathcal{L}^{\mathrm{4D}}_{2})$ and $U(\mathcal{L}^{\mathrm{4D}}_{3})$ 
are as follows:
\begin{description}
\item[$U(\mathcal{L}^{\mathrm{4D}}_{1})$]
\begin{align}
\chi_{3}\chi_{4}&=\frac{z(h+g)}{z+1}\chi_{1}^{2}+(h+g)^{2}\chi_{2}\chi_{1}-
(h+g)\chi_{2}\chi_{3}+\frac{z(h+g)}{z+1}\chi_{1}\chi_{4}-
\frac{z}{z+1}\chi_{1}\chi_{3}\nonumber\\
&+\frac{1}{z+1}\chi_{4}\chi_{3}+
\frac{h+g}{z+1}\chi_{1}-\frac{1}{z+1}\chi_{3},\nonumber\\
\chi_{3}\chi_{1}&=-\frac{z(h+g)}{z+1}\chi_{1}^{2}-(h+g)^{2}\chi_{2}\chi_{1}+
(h+g)\chi_{2}\chi_{3}-\frac{z(h+g)}{z+1}\chi_{1}\chi_{4}+
\frac{2z+1}{z+1}\chi_{1}\chi_{3}\nonumber\\
&+\frac{z}{z+1}\chi_{4}\chi_{3}-
\frac{h+g}{z+1}\chi_{1}+\frac{1}{z+1}\chi_{3},\nonumber\\
\chi_{3}\chi_{2}&=\frac{z}{z+1}\chi_{4}^{2}-\frac{z}{z+1}\chi_{1}^{2}-
(h+g)^{2}\chi_{2}^{2}-\frac{(2z+1)(h+g)}{z+1}\chi_{2}\chi_{1}+
\frac{h+g}{z+1}\chi_{2}\chi_{4}+\chi_{2}\chi_{3}\nonumber\\
&-\frac{1}{z+1}\chi_{1}-
\frac{h+g}{z+1}\chi_{2}+\frac{1}{z+1}\chi_{4},\nonumber\\
\chi_{4}\chi_{1}&=\chi_{1}\chi_{4},\nonumber\\
\chi_{4}\chi_{2}&=-(h+g)\chi_{2}^{2}-\frac{z}{z+1}\chi_{2}\chi_{1}+
\frac{1}{z+1}\chi_{2}\chi_{4}-\frac{1}{z+1}\chi_{2},\nonumber\\
\chi_{1}\chi_{2}&=(h+g)\chi_{2}^{2}+\frac{2z+1}{z+1}\chi_{2}\chi_{1}+
\frac{z}{z+1}\chi_{2}\chi_{4}+\frac{1}{z+1}\chi_{2}.
\end{align}
\item[$U(\mathcal{L}^{\mathrm{4D}}_{2})$]
\begin{align}
\chi_{3}\chi_{4}&=(h+g)^{2}\chi_{2}\chi_{1}-(h+g)\chi_{2}\chi_{3}+
\chi_{4}\chi_{3}+(h+g)\chi_{1}-\chi_{3},\nonumber\\
\chi_{3}\chi_{1}&=-(h+g)^{2}\chi_{2}\chi_{1}+(h+g)\chi_{2}\chi_{3}+
\chi_{1}\chi_{3}-(h+g)\chi_{1}+\chi_{3},\nonumber\\
\chi_{3}\chi_{2}&=-(h+g)^{2}\chi_{2}^{2}-(h+g)\chi_{2}\chi_{1}+
(h+g)\chi_{2}\chi_{4}+\chi_{2}\chi_{3}-\chi_{1}-(h+g)\chi_{2}+\chi_{4},\nonumber\\
\chi_{4}\chi_{1}&=\chi_{1}\chi_{4},\nonumber\\
\chi_{4}\chi_{2}&=-(h+g)\chi_{2}^{2}+\chi_{2}\chi_{4}-\chi_{2},\nonumber\\
\chi_{1}\chi_{2}&=(h+g)\chi_{2}^{2}+\chi_{2}\chi_{1}+\chi_{2}.
\end{align}
\item[$U(\mathcal{L}^{\mathrm{4D}}_{3})$]
\begin{align}
\chi_{3}\chi_{4}&=(h+g)^{2}\chi_{2}\chi_{1}-(h+g)\chi_{2}\chi_{3}+
\chi_{4}\chi_{3}+(h+g)\chi_{1}-\chi_{3},\nonumber\\
\chi_{3}\chi_{1}&=-(h+g)^{2}\chi_{2}\chi_{1}+(h+g)\chi_{2}\chi_{3}+
\chi_{1}\chi_{3}-(h+g)\chi_{1}+\chi_{3},\nonumber\\
\chi_{3}\chi_{2}&=-(h+g)^{2}\chi_{2}^{2}-(h+g)\chi_{2}\chi_{1}+
(h+g)\chi_{2}\chi_{4}+\chi_{2}\chi_{3}-\chi_{1}-(h+g)\chi_{2}+\chi_{4},\nonumber\\
\chi_{4}\chi_{1}&=\chi_{1}\chi_{4},\nonumber\\
\chi_{4}\chi_{2}&=-(h+g)\chi_{2}^{2}+\chi_{2}\chi_{4}-\chi_{2},\nonumber\\
\chi_{1}\chi_{2}&=(h+g)\chi_{2}^{2}+\chi_{2}\chi_{1}+\chi_{2}.
\end{align}
\end{description}
In each case the relations are such that 
$\{\chi_{2}^{\a}\chi_{1}^{\b}\chi_{4}^{\g}\chi_{2}^{\d}:\a,\b,\g,\d\in\ZZ_{\geq 0}\}$
is a basis of the enveloping algebra.
\ethe
\begin{proof}
These relations are obtained by solving the 16 equations $\chi_{i}\chi_{k}=\sum_{s,t=1\ldots 
d}\L_{st,ik}\chi_{s}\chi_{t}+\sum_{j=1\ldots 
d}\mathcal{C}_{ik,j}\chi_{j}$ for the 6 quadratic elements 
$\chi_{3}\chi_{4}$, $\chi_{3}\chi_{1}$, $\chi_{3}\chi_{2}$, 
$\chi_{4}\chi_{1}$, $\chi_{4}\chi_{2}$ and $\chi_{1}\chi_{2}$. It is 
then observed that the relations are compatible with the ordering 
$\chi_{2}\prec\chi_{1}\prec\chi_{4}\prec\chi_{3}$ so the Diamond Lemma 
may be applied to obtain the stated basis.
\end{proof}

\label{ClasResults1}

\section{3-dimensional bicovariant differential calculi on $SL_{h}(2)$}
Classically (see for example the discussion in the book by 
Flanders~\cite{Flanders}) we obtain the differential calculus on 
$SL(2)$ from the calculus on $GL(2)$ through the classical relation,
\begin{equation}
\dif\mathcal{D}=\mathcal{D}\Tr\Theta^{c},\label{cldifdet}
\end{equation}
where $\Theta^{c}$ is the classical matrix of left-invariant 1-forms.
$SL(2)$ and its differential calculus is obtained by setting $\mathcal{D}=1$,
so the left hand side 
becomes zero and we obtain a linear relation between the classical
left-invariant 1-forms, namely, $\t_{1}+\t_{4}=0$. In the standard 
quantum $GL_{q}(2)$ case this procedure is \emph{not} possible. The analogue 
of~(\ref{cldifdet}) in this case is rendered trivial by the condition that 
$\mathcal{D}$ be central in the first order calculus. More precisely, 
we have in this case
\begin{equation}
\dif\mathcal{D}=\k\mathcal{D}\Tr_{q}\Theta,
\end{equation}
where $\Tr_{q}\Theta$ is now the $q$-analogue of $\Tr\Theta$, 
$\t_{1}+q^{-1}\t_{4}$. But now, imposing the condition $\mathcal{D}=1$, 
immediately fixes $\k=0$ so we have no chance of reducing the 
dimension of the calculus. There are of course 4-dimensional 
bicovariant calculi 
on the standard quantum group $SL_{q}(2)$ but \emph{no} 3-dimensional 
calculi.

In our case, with the non-standard Jordanian quantum group, the 
situation is quite different. Studying Theorem~\ref{glabcd} we see 
that for the calculi $\G^{\mathrm{4D}}_{1}$ and $\G^{\mathrm{4D}}_{3}$ 
we can indeed have the quantum determinant central \emph{and} obtain a 
dimension reducing relation through the analogue in these cases 
of~(\ref{cldifdet}),
\begin{equation}
\dif\mathcal{D}=\frac{z+2}{2}\mathcal{D}\Tr_{h}\Theta.
\end{equation}
The quantum determinant is central in $\G^{\mathrm{4D}}_{1}$ 
and $\G^{\mathrm{4D}}_{3}$ if and only if the parameter $z$ takes the 
value $0$ or $-2$. With $z=-2$ we obtain in each case a 4-dimensional 
calculus on $SL_{h}(2)$, while with $z=0$ the condition 
$\dif\mathcal{D}=0$ yields the linear relation 
$\Tr_{h}\Theta=\t_{1}+2h\t_{3}+\t_{4}=0$ --- precisely the 
relation~(\ref{3Dred}) we 
obtained when we investigated the general implications of choosing a 
three dimensional basis of $\G_{\inv}$. Moreover, $z=0$ is the value 
of $z$ at which $\G^{\mathrm{4D}}_{1}$, $\G^{\mathrm{4D}}_{2}$ and 
$\G^{\mathrm{4D}}_{3}$ coincide, so is already covered as a 
particular case in 
$\G^{\mathrm{4D}}_{2}$. For this first order calculus, having set 
$g=h$, the quantum determinant is central for \emph{all} values of z. 
In the particular case of $z=\frac{2}{3}$, the differential of the 
quantum determinant is identically zero so that once again we obtain a 
4-dimensional calculus on $SL_{h}(2)$. However, for all other values 
of $z$ we recover the condition $\Tr_{h}\Theta=0$. At first sight then, it may 
seem that there is a family of 3-dimensional calculi on 
$SL_{h}(2)$, but this is not the case.

\bth
There is a unique, 3-dimensional, first order bicovariant 
differential calculus on the Jordanian quantum group $SL_{h}(2)$, 
$\G^{\mathrm{3D}}$. It 
may be obtained from any one of the three families of first order bicovariant 
differential calculi on $GL_{h,g}(2)$ by a reduction analogous to the 
classical situation. It is specified by its $ABCD$ matrices,
\begin{align}
A&=\begin{pmatrix}
1&0&-h\\
2h&1&h^{2}\\
0&0&1
\end{pmatrix},&
B&=\begin{pmatrix}
0&0&2h^{2}\\
0&2h&-2h^{3}\\
0&0&-2h
\end{pmatrix}\nonumber\\
C&=\begin{pmatrix}
0&0&0\\
0&0&0\\
0&0&0
\end{pmatrix},&
D&=\begin{pmatrix}
1&0&h\\
-2h&1&-3h^{2}\\
0&0&1
\end{pmatrix}.
\end{align}
\ethe
\begin{proof}
As far as obtaining this calculus from the $GL_{h,g}(2)$ calculi is 
concerned, we 
need only observe that starting with the $4\times 4$ $ABCD$ matrices of 
$\G^{\mathrm{4D}}_{3}$, say, in the relations~(\ref{thfmcrs}) and 
setting $\t_{4}=-\t_{1}-2h\t_{3}$, we obtain commutation relations now 
involving the $3\times 3$ matrices quoted. That this is the unique 
3-dimensional calculus on $SL_{h}(2)$ follows since these are 
precisely the $ABCD$ matrices we obtain when we apply the 
procedure of Section~4 to look for the most general possible 
3-dimensional calculus on $SL_{h}(2)$.
\end{proof}

Now, just as we did for the 4D calculi, we may deduce wedge product 
commutation relations in the exterior bicovariant graded algebra
$\Omega^{\mathrm{3D}}$, the Lie brackets and commutators for the 
quantum Lie algebra $\mathcal{L}^{\mathrm{3D}}$, and the 
enveloping algebra relations for $U(\mathcal{L}^{\mathrm{3D}})$. These 
results are obtained in just the same way as the corresponding results 
for the 4D calculi so we will not comment on their proofs. We 
should mention that some results in this direction have already been 
obtained in~\cite{Kar2}, where the first order differential calculus 
on $SL_{h}(2)$ was postulated through the $R$-matrix 
expression~(\ref{rmatdif}). However the results we will present here 
are more complete than the corresponding results 
in~\cite{Kar2}. Also, it will be useful to collect the results here 
in our current notation.

\bth
The commutation relations between the left-invariant forms in the 
external bicovariant graded algebra, $\Omega^{\mathrm{3D}}$, built 
upon the first order calculus $\G^{\mathrm{3D}}$, are
\begin{align}
\t_{3}\wedge\t_{3}&=0,\nonumber\\
\t_{3}\wedge\t_{1}&=-\t_{1}\wedge\t_{3},\nonumber\\
\t_{3}\wedge\t_{2}&=-\t_{2}\wedge\t_{3}+4h\t_{1}\wedge\t_{3},\nonumber\\
\t_{1}\wedge\t_{1}&=0,\nonumber\\
\t_{1}\wedge\t_{2}&=-\t_{2}\wedge\t_{1}-2h\t_{2}\wedge\t_{3},\nonumber\\
\t_{2}\wedge\t_{2}&=4h\t_{2}\wedge\t_{1}+8h^{2}\t_{2}\wedge\t_{3}.
\end{align}
The relations are such that 
$\{\t_{2}^{\a}\t_{1}^{\b}\t_{3}^{\g}:\a,\b,\g\in\{0,1\}\}$ is a basis 
for the exterior algebra of forms.
\ethe

\bth
The quantum Lie brackets and commutators for the quantum Lie algebra 
$\mathcal{L}^{\mathrm{3D}}$, are
respectively,
\begin{align*}
[\chi_{1},\chi_{1}]&=0,\nonumber\\
[\chi_{1},\chi_{2}]&=2\chi_{2},\nonumber\\
[\chi_{1},\chi_{3}]&=-2\chi_{3},\nonumber\\
[\chi_{2},\chi_{1}]&=-2\chi_{2},\nonumber\\
[\chi_{2},\chi_{2}]&=0,\nonumber\\
[\chi_{2},\chi_{3}]&=\chi_{1}-4h\chi_{2},\nonumber\\
\end{align*}
\begin{align}
[\chi_{3},\chi_{1}]&=-4h\chi_{1}+2\chi_{3},\nonumber\\
[\chi_{3},\chi_{2}]&=-\chi_{1},\nonumber\\
[\chi_{3},\chi_{3}]&=-4h\chi_{3},\label{jqlab}
\end{align}
and 
\begin{align}
[\chi_{1},\chi_{1}]&=0,\nonumber\\
[\chi_{1},\chi_{2}]&=\chi_{1}\chi_{2}-\Bigl(\chi_{2}\chi_{1}+
4h\chi_{2}\chi_{2}\Bigr),\nonumber\\
[\chi_{1},\chi_{3}]&=\chi_{1}\chi_{3}-\Bigl(\chi_{3}\chi_{1}-
4h\chi_{3}\chi_{2}\Bigr),\nonumber\\
[\chi_{2},\chi_{1}]&=\chi_{2}\chi_{1}-\Bigl(\chi_{1}\chi_{2}-
4h\chi_{2}\chi_{2}\Bigr),\nonumber\\
[\chi_{2},\chi_{2}]&=0,\nonumber\\
[\chi_{2},\chi_{3}]&=\chi_{2}\chi_{3}-\Bigl(2h\chi_{1}\chi_{2}-
8h^{2}\chi_{2}\chi_{2} + \chi_{3}\chi_{2}\Bigr),\nonumber\\
[\chi_{3},\chi_{1}]&=\chi_{3}\chi_{1}-\Bigl(\chi_{1}\chi_{3}-
8h^{2}\chi_{2}\chi_{1}+4h\chi_{2}\chi_{3}\Bigr),\nonumber\\
[\chi_{3},\chi_{2}]&=\chi_{3}\chi_{2}-\Bigl(-2h\chi_{2}\chi_{1}+
\chi_{2}\chi_{3}\Bigr),\nonumber\\
[\chi_{3},\chi_{3}]&=\chi_{3}\chi_{3}-\Bigl(-2h\chi_{1}\chi_{3}+
2h\chi_{3}\chi_{1}-8h^{2}\chi_{3}\chi_{2}+\chi_{3}\chi_{3}\Bigr).
\end{align}
\ethe

\bth
The relations in the enveloping algebra, 
$U(\mathcal{L}^{\mathrm{3D}})$, are,
\begin{align}
\chi_{3}\chi_{1}&=-8h^{2}\chi_{2}\chi_{1}+4h\chi_{2}\chi_{3}+\chi_{1}\chi_{3}-
4h\chi_{1}+2\chi_{3},\nonumber\\
\chi_{3}\chi_{2}&=-2h\chi_{2}\chi_{1}+\chi_{2}\chi_{3}-\chi_{1},\nonumber\\
\chi_{1}\chi_{2}&=4h\chi_{2}^{2}+\chi_{2}\chi_{1}+2\chi_{2}.
\end{align}
These relations are such that 
$\{\chi_{2}^{\a}\chi_{1}^{\b}\chi_{3}^{\g}:\a,\b,\g\in\ZZ_{\geq 0}\}$ 
is a basis for $U(\mathcal{L}^{\mathrm{3D}})$.
\ethe

\label{ClasResults2}

\section{The Jordanian quantised universal enveloping algebra}
To this point we have been working only with the Jordanian quantum 
analog of the coordinate ring of $SL_{2}(\CC)$. In the remainder of 
the paper we focus attention on the corresponding deformation of the 
universal enveloping algebra $U(\slt)$. Let us recall its
definition.

\bde
The \emph{Jordanian quantised universal enveloping algebra,} 
$U_{h}(\slt)$, is the unital associative algebra over $\CC[[h]]$ with 
generators $X$, $Y$, $H$ and relations
\begin{gather}
[H,X]=2\frac{\sinh hX}{h},\\
[H,Y]=-Y(\cosh hX)-(\cosh hx)Y,\\
[X,Y]=H.
\end{gather}
having a basis $\{Y^{\a}H^{\b}X^{\g}:\a,\b,\g\in\ZZ_{\geq 0}\}$.
\ede

In~\cite{Ang} the Casimir element, $C$, of $U_{h}(\slt)$ was obtained,
in our notation it is
\begin{equation}
C=\left(Y(\sinh hX)+(\sinh 
hX)Y\right)+\frac{1}{4}H^{2}+\frac{1}{4}(\sinh hX)^{2}.
\end{equation}

Assuming tensor products to be completed, where appropriate, in the 
$h$-adic topology, the Hopf structure of $U_{h}(\slt)$ is defined on 
the generators as,
\begin{gather}
\D(X)=X\ot 1+1\ot X,\\
\D(Y)=Y\ot e^{hX}+e^{-hX}\ot Y,\\
\D(H)=H\ot e^{hX}+e^{-hX}\ot H,\\
\ep(X)=0,\;\;\;\;\ep(Y)=0,\;\;\;\;\ep(H)=0,\\
S(X)=-X,\;\;\;\;S(Y)=-e^{hX}Ye^{-hX},\;\;\;\;S(H)=-e^{hX}He^{-hX}.
\end{gather}

$U_{h}(\slt)$ has received a good deal of 
attention recently. In particular, we mention the work of 
Abdesselam~\emph{et al}~\cite{ACC} in which a non-linear map was constructed 
which realises the \emph{algebraic} isomorphism between $U_{h}(\slt)$ and 
$U(\slt)$. This map was then used by those authors to build the representation 
theory of $U_{h}(\slt)$. As with the standard quantisation of the 
enveloping algebra of $\slt$, the representation theory of 
$U_{h}(\slt)$ follows very closely the representation theory of $\slt$. 
Indeed, the finite dimensional, indecomposable representations of 
$U_{h}(\slt)$ are in one-to-one correpondence with the finite 
dimensional irreducible representations of $\slt$, and can be 
classified as classically by a non-negative half-integer $j$. Van der 
Jeugt~\cite{vdJ} was able to refine the work of Abdesselam~\emph{et al}, 
obtaining closed form expressions for the action of the 
generators of $U_{h}(\slt)$ on the basis vectors of finite 
dimensional irreducible representations. Before Van der 
Jeugt's work Aizawa~\cite{Aiz} had demonstrated that the Clebsch-Gordan series 
for the decomposition of the tensor product of two indecomposable representations
of $U_{h}(\slt)$ was 
precisely the classical series modulo the one-to-one correspondence 
of classical and Jordanian representations. Van der Jeugt obtained a 
general formula for the Clebsch-Gordan coefficients.

\label{defuhslt}

\section{Jordanian quantum Lie algebra from an ad-submodule in 
$U_{h}(\slt)$}

In the following theorem we describe a left ad-submodule of 
$U_{h}(\slt)$ which allows us to build a quantum Lie algebra from the 
enveloping algebra generators.
\bth
In $U_{h}(\slt)$ the space spanned by the elements $X_{h}$, $H_{h}$ and 
$Y_{h}$ defined by
\begin{align}
X_{h}&=e^{hX}\frac{\sinh hX}{h},\nonumber\\
H_{h}&=He^{hX},\nonumber\\
Y_{h}&=Ye^{hX}-2hC,
\end{align}
is stable under the left adjoint action of $U_{h}(\slt)$ on $U_{h}(\slt)$.
\ethe
\begin{proof}
To obtain this result, essential use was made of the known PBW basis. 
With such a basis we can use the computer algebra package REDUCE to 
perform algebraic manipulations which would be virtually impossible 
otherwise. 
\end{proof}

The actions of the $U_{h}(\slt)$ generators on the elements 
$\{X_{h},H_{h},Y_{h}\}$ is given 
by,
\begin{align}
X\act X_{h}&=0,&H\act X_{h}&=2X_{h},&Y\act 
X_{h}&=-H_{h}+2hX_{h},\nonumber\\
X\act H_{h}&=-2X_{h},&H\act H_{h}&=4hX_{h},&Y\act 
H_{h}&=2Y_{h}+3h^{2}X_{h},\nonumber\\
X\act Y_{h}&=H_{h},&H\act Y_{h}&=-2Y_{h}-2hH_{h}-h^{2}X_{h},
&Y\act Y_{h}&=-2hY_{h}-h^{2}H_{h}-h^{3}X_{h}\label{rep1}
\end{align}
The actions of the elements on each other leads to the following 
\emph{Jordanian quantum Lie brackets},
\begin{align}
[X_{h},X_{h}]&=0,&[X_{h},H_{h}]&=-2X_{h},&[X_{h},Y_{h}]&=H_{h}-2hX_{h},\nonumber\\
[H_{h},X_{h}]&=2X_{h},&[H_{h},H_{h}]&=0,&[H_{h},Y_{h}]&=-2Y_{h}-2hH_{h}+h^{2}X_{h},\nonumber\\
[Y_{h},X_{h}]&=-H_{h}-2hX_{h},&[Y_{h},H_{h}]&=2Y_{h}-2hH_{h}-h^{2}X_{h},
&[Y_{h},Y_{h}]&=-4hY_{h},\label{Jqlbr}
\end{align}
which display the characteristic $h$-antisymmetry~\cite{Del2}.
The $U_{h}(\slt)$ coproduct on the elements $X_{h}$, $H_{h}$ and $Y_{h}$ 
is,
\begin{align}
\D(X_{h})&=1\ot X_{h}+X_{h}\ot e^{2hX},\nonumber\\
\D(H_{h})&=1\ot H_{h}+H_{h}\ot e^{2hX},\nonumber\\
\D(Y_{h})&=1\ot Y_{h}+Y_{h}\ot e^{2hX} +2h(1\ot C + C\ot e^{2hX} - 
\D(C)).
\end{align}
\label{Approach1}

\section{Jordanian quantum Lie algebra from inverse Clebsch-Gordan 
coefficients}
In the notation of Van der Jeugt~\cite{vdJ}, on the representation 
space $V^{(j)}_{h}$ with basis $e^{j}_{m}$ where $j$ in a non-negative 
half-integer and $m=-j,-j+1,\ldots,j$, the action of the generators 
$X$, $H$ and $Y$ is given by~\cite{vdJ},
\begin{align}
H\act e^{j}_{m}&=2me^{j}_{m},\nonumber\\
X\act e^{j}_{m}&=
\sum_{k=0}^{\lfloor(j-m-1)/2\rfloor}
\frac{(h/2)^{2k}}{2k+1}\frac{\a_{j,m+1+2k}}{\a_{j,m}}e^{j}_{m+1+2k},\nonumber\\
Y\act 
e^{j}_{m}&=(j+m)(j-m+1)\frac{\a_{j,m-1}}{\a_{j,m}}e^{j}_{m-1}\nonumber\\
&\quad-
(j-m)(j+m+1)\left(\frac{h}{2}\right)^{2}\frac{\a_{j,m+1}}{\a_{j,m}}e^{j}_{m+1}
\nonumber\\
&\quad
+\sum_{s=1}^{\lfloor(j-m+1)/2\rfloor}\left(\frac{h}{2}\right)^{2s}
\frac{\a_{j,m-1+2s}}{\a_{j,m}}e^{j}_{m-1+2s},
\end{align}
where $\a_{j,m}=\sqrt{(j+m)!/(j-m)!}$.
Thus the representation matrices of the generators in the deformation of 
the classical adjoint, $j=1$, representation are,
\begin{gather}
\Gamma(X)=\begin{pmatrix}
0&\sqrt{2}&0\\
0&0&\sqrt{2}\\
0&0&0
\end{pmatrix},\;\;\;\;
\Gamma(H)=\begin{pmatrix}
2&0&0\\
0&0&0\\
0&0&-2
\end{pmatrix},\nonumber\\
\Gamma(Y)=\begin{pmatrix}
0&-\sqrt{2}(h/2)^{2}&0\\
\sqrt{2}&0&-\sqrt{2}(h/2)^{2}\\
0&\sqrt{2}&0
\end{pmatrix}.
\end{gather}

The Clebsch-Gordan series for the tensor product of two $V^{1}_{h}$ 
representations is,
\begin{equation}
V^{1}_{h}\ot V^{1}_{h}\cong V^{2}_{h}\opl V^{1}_{h}\opl V^{0}_{h}.
\end{equation}
If we denote by $v_{i}$ the $i$-th vector in the ordered basis 
$\{e^{2}_{2},e^{2}_{1},e^{2}_{0},e^{2}_{-1},e^{2}_{-2},e^{1}_{1},e^{1}_{0},
e^{1}_{-1},e^{0}_{0}\}$ and by $w_{i}$ the $i$-th vector in the 
ordered basis $\{e^{1}_{1}\ot e^{1}_{1},e^{1}_{1}\ot e^{1}_{0},
e^{1}_{1}\ot e^{1}_{-1},e^{1}_{0}\ot e^{1}_{1},e^{1}_{0}\ot e^{1}_{0},
e^{1}_{0}\ot e^{1}_{-1},e^{1}_{-1}\ot e^{1}_{1},e^{1}_{-1}\ot e^{1}_{0},
e^{1}_{-1}\ot e^{1}_{-1}\}$, then the Clebsch-Gordan matrix $C$, where
$v_{i}=\sum_{j=1}^{9}C_{ij}w_{j}$ is given by~\cite{vdJ}
\begin{equation}
C=\begin{pmatrix}
1 & 0 & 0 & 0 & 0 & 0 & 0 & 0 & 0 \\
0 & \frac{1}{\sqrt{2}} & 0 & \frac{1}{\sqrt{2}} & 0 & 0 & 0 & 0 & 0 \\
\frac{\sqrt{2} h^{2}}{2 \sqrt{3}} & \frac{-h}{\sqrt{3}} & \frac{1}{\sqrt{6}} & 
\frac{h}{\sqrt{3}} & \frac{2}{\sqrt{6}} & 0 & \frac{1}{\sqrt{6}} & 0 & 0 \\
0 & \frac{\sqrt{2} h^{2}}{2} & -h & \frac{\sqrt{2} h^{2}}{2} & 0 & \frac{1}{\sqrt{2}} & h & 
\frac{1}{\sqrt{2}} & 0\\
\frac{-h^{4}}{4} & \frac{-\sqrt{2} h^{3}}{2} & \frac{3h^{2}}{2} & 
\frac{\sqrt{2} h^{3}}{2} & 0 & -\sqrt{2} h & \frac{3h^{2}}{2} & \sqrt{2} h & 1 \\
-2h & \frac{1}{\sqrt{2}} & 0 & \frac{-1}{\sqrt{2}} & 0 & 0 & 0 & 0 & 
0 \\
0 & -h & \frac{1}{\sqrt{2}} & -h & 0 & 0 & \frac{-1}{\sqrt{2}} & 0 & 0 \\
0 & \frac{\sqrt{2}h^{2}}{2} & -h & \frac{-\sqrt{2} h^{2}}{2} & 0 & \frac{1}{\sqrt{2}} 
& -h & \frac{-1}{\sqrt{2}} & 0 \\
\frac{h^{2}}{\sqrt{3}} & \frac{-\sqrt{2}h}{\sqrt{3}} & 
\frac{1}{\sqrt{3}} & \frac{\sqrt{2} h}{\sqrt{3}} & \frac{-1}{\sqrt{3}} & 
0 & \frac{1}{\sqrt{3}} & 0 & 0
\end{pmatrix},
\end{equation}
with inverse
\begin{equation}
C^{-1}=\begin{pmatrix}
1 & 0 & 0 & 0 & 0 & 0 & 0 & 0 & 0 \\
\sqrt{2}h & \frac{\sqrt{2}}{2} & 0 & 0 & 0 & \frac{\sqrt{2}}{2} & 0 & 0 & 0 \\
\frac{3h^{2}}{2} & h & \frac{\sqrt{6}}{6} & 0 & 0 & h & 
\frac{\sqrt{2}}{2} & 0 & \frac{\sqrt{3}}{3} \\
-\sqrt{2}h & \frac{\sqrt{2}}{2} & 0 & 0 & 0 & \frac{-\sqrt{2}}{2} & 0 & 0 & 0 \\
0 & 0 & \frac{\sqrt{6}}{3} & 0 & 0 & 0 & 0 & 0 & \frac{-\sqrt{3}}{3} \\
\frac{\sqrt{2}h^{3}}{2} & \frac{\sqrt{2}h^{2}}{2} & \frac{\sqrt{3}h}{3} & 
\frac{\sqrt{2}}{2} & 0 & \frac{\sqrt{2}h^{2}}{2} & h & 
\frac{\sqrt{2}}{2} & \frac{\sqrt{6}h}{3} \\
\frac{3h^{2}}{2} & -h & \frac{\sqrt{6}}{6} & 0 & 0 & h & 
\frac{-\sqrt{2}}{2} & 0 & \frac{\sqrt{3}}{3} \\
\frac{-\sqrt{2}h^{3}}{2} & \frac{\sqrt{2}h^{2}}{2} & 
\frac{-\sqrt{3}h}{3} & \frac{\sqrt{2}}{2} & 0 & 
\frac{-\sqrt{2}h^{2}}{2} & h & \frac{-\sqrt{2}}{2} & 
\frac{-\sqrt{6}h}{3} \\
\frac{-h^{4}}{4} & 0 & \frac{\sqrt{6}h^{2}}{6} & 0 & 1 & 0 & 0 & 
2h & \frac{\sqrt{3}h^{2}}{3}
\end{pmatrix}.
\end{equation}

Considering $w_{i}=\sum_{j=1}^{9}C^{-1}_{ij}v_{j}$ we see that columns 6--8 of 
$C^{-1}$ correspond to an intertwiner 
$\mathrm{ad}\ot\mathrm{ad}\mapto\mathrm{ad}$ and 
we deduce a quantum Lie bracket on the 
vectors $\{e^{1}_{1},e^{1}_{0},e^{1}_{-1}\}$,
\begin{align}
[e^{1}_{1},e^{1}_{1}]&=0,&[e^{1}_{1},e^{1}_{0}]&=\frac{\sqrt{2}}{2}e^{1}_{1},&
[e^{1}_{1},e^{1}_{-1}]&=\frac{\sqrt{2}}{2}e^{1}_{0}+he^{1}_{1},\nonumber\\
[e^{1}_{0},e^{1}_{1}]&=-\frac{\sqrt{2}}{2}e^{1}_{1},&
[e^{1}_{0},e^{1}_{0}]&=0,&
[e^{1}_{0},e^{1}_{-1}]&=\frac{\sqrt{2}}{2}e^{1}_{-1}+he^{1}_{0}+
\frac{\sqrt{2}}{2}h^{2}e^{1}_{1},\nonumber\\
[e^{1}_{-1},e^{1}_{1}]&=-\frac{\sqrt{2}}{2}e^{1}_{0}+he^{1}_{1},&
[e^{1}_{-1},e^{1}_{0}]&=-\frac{\sqrt{2}}{2}e^{1}_{-1}+he^{1}_{0}-
\frac{\sqrt{2}}{2}h^{2}e^{1}_{1},&
[e^{1}_{-1},e^{1}_{-1}]&=2he^{1}_{-1}.
\end{align}
Similarly, column 9 of $C^{-1}$ corresponds to an intertwiner 
$\mathrm{ad}\ot\mathrm{ad}\mapto\CC[[h]]$ and we obtain the Killing form
\begin{align}
\kf(e^{1}_{1},e^{1}_{-1})&=\frac{\sqrt{3}}{3},&
\kf(e^{1}_{-1},e^{1}_{1})&=\frac{\sqrt{3}}{3},\nonumber\\
\kf(e^{1}_{0},e^{1}_{0})&=-\frac{\sqrt{3}}{3},&
\kf(e^{1}_{-1},e^{1}_{0})&=-\frac{\sqrt{6}h}{3},\nonumber\\
\kf(e^{1}_{0},e^{1}_{-1})&=\frac{\sqrt{6}h}{3},&
\kf(e^{1}_{-1},e^{1}_{-1})&=\frac{\sqrt{3}h^{2}}{3},
\end{align}

Now, if we perform the following change of basis,
\begin{align}
X_{h}&=2e^{1}_{1},\nonumber\\
H_{h}&=4he^{1}_{1}-2\sqrt{2}e^{1}_{0},\nonumber\\
Y_{h}&=-\frac{5}{2}h^2e^{1}_{1}+2\sqrt{2}he^{1}_{0}-2e^{1}_{-1},
\end{align}
then the representation matrices of the generators become,
\begin{gather}
\Gamma(X)=\begin{pmatrix}
0&-2&0\\
0&0&1\\
0&0&0
\end{pmatrix},\;\;\;\;
\Gamma(H)=\begin{pmatrix}
2&4h&-h^{2}\\
0&0&-2h\\
0&0&-2
\end{pmatrix},\nonumber\\
\Gamma(Y)=\begin{pmatrix}
2h&3h^{2}&-h^{3}\\
-1&0&-h^{2}\\
0&2&-2h
\end{pmatrix},\label{adreps}
\end{gather}
which are precisely those obtained from~(\ref{rep1}) 
and the Lie bracket relations become (those already obtained 
in~\ref{Jqlbr},
\begin{align}
[X_{h},X_{h}]&=0,&[X_{h},H_{h}]&=-2X_{h},&[X_{h},Y_{h}]&=H_{h}-2hX_{h},\nonumber\\
[H_{h},X_{h}]&=2X_{h},&[H_{h},H_{h}]&=0,&[H_{h},Y_{h}]&=-2Y_{h}-2hH_{h}+h^{2}X_{h},\nonumber\\
[Y_{h},X_{h}]&=-H_{h}-2hX_{h},&[Y_{h},H_{h}]&=2Y_{h}-2hH_{h}-h^{2}X_{h},
&[Y_{h},Y_{h}]&=-4hY_{h}.
\end{align}
Further, we can scale the Killing form, by scaling the single basis 
vector of $V^{0}_{h}$ so that on the $\{X_{h},H_{h},Y_{h}\}$ it reads,
\begin{align}
\kf(X_{h},Y_{h})&=4,&
\kf(Y_{h},X_{h})&=4,\nonumber\\
\kf(H_{h},H_{h})&=8,&
\kf(Y^{h},H_{h})&=-8h,\nonumber\\
\kf(H_{h},Y_{h})&=-8h,&
\kf(Y_{h},Y_{h})&=-6h^{2},
\end{align}
a simple deformation of the classical Killing form, recovered by 
setting $h=0$.
\label{Approach2}

\section{Conclusion}
Returning to the quantum Lie algebra obtained through Woronowicz's 
bicovariant calculus, $\mathcal{L}^{\mathrm{3D}}$, 
if we change the basis according to the 
identifications,
\begin{align}
H_{h}&=\chi_{1},\nonumber\\
X_{h}&=\chi_{2},\nonumber\\
Y_{h}&=-h\chi_{1}+\frac{h^{2}}{4}\chi_{2}+\chi_{3}
\end{align}
then the \emph{Woronowicz quantum Lie bracket} on these new basis elements 
is \emph{precisely} that already found in~(\ref{Jqlbr}). Thus, as 
algebras over $\CC[[h]]$, the Woronowicz and `Sudbery-Delius' quantum 
Lie algebras are \emph{isomorphic}. This means, furthermore, that in 
addition to having a Killing form we have some analog of the Jacobi 
identity for this Jordanian quantum lie algebra.

We had already found two appealing aspects of the bicovariant 
differential geometry on $SL_{h}(2)$. Namely, its uniqueness and 
3-dimensionality. The fact that the Woronowicz quantum Lie algebra is 
isomorphic to the Sudbery-Delius quantum Lie algebra we found 
starting with $U_{h}(\slt)$ is a further attractive feature. 
Recently, Cho,~Madore 
and~Park~\cite{CMP} have shown that the corresponding Jordanian 
quantum plane 
admits a richer geometrical structure than the standard quantum 
plane. It should be interesting to try to develop further the geometry on the 
Jordanian quantum group and also investigate possible $SL(n)$ 
generalisations.
\label{last}

\end{document}